%% file: jinv.tex
\documentclass [10pt,reqno]{amsart}
\usepackage {hannah}
\usepackage {color,graphicx}
\usepackage {multirow}
\usepackage{texdraw}
\usepackage [latin1]{inputenc}

\newcommand {\RR}{{\mathbb R}}
\newcommand {\ZZ}{{\mathbb Z}}

\DeclareMathOperator {\ev}{ev}

\DeclareMathOperator {\val}{val}

\DeclareMathOperator {\defi}{def}
\DeclareMathOperator {\trop}{trop}

\DeclareMathOperator {\Area}{Area}

\newcommand {\dunion}{\,\mbox {\raisebox{0.25ex}{$\cdot$} \kern-1.83ex $\cup$}
  \,}

\title [Counting tropical elliptic plane curves with fixed $j$-invariant]{Counting tropical elliptic plane curves with fixed $j$-invariant}
\author {Michael Kerber and Hannah Markwig}
\address {Michael Kerber, Fachbereich Mathematik, Technische Universit\"at
  Kaiserslautern, Postfach 3049, 67653 Kaiserslautern, Germany}
\email {mkerber@mathematik.uni-kl.de}
\address {Hannah Markwig, Fachbereich Mathematik, Technische Universit\"at Kaiserslautern, Postfach 3049, 67653 Kaiserslautern, Germany}
\email {markwig@mathematik.uni-kl.de}

\begin {document}
\begin {abstract}
   In complex algebraic geometry, the problem of enumerating plane
   elliptic curves of given degree with fixed complex structure has been
   solved by R.Pandharipande \cite{Pan97} using
   Gromov-Witten theory. In this article we treat the tropical
   analogue of this problem, the determination of the number $E_{\trop}(d)$ of
   tropical elliptic plane curves of degree $d$ and fixed ``tropical
   $j$-invariant'' interpolating an appropriate number of points in general position and counted with multiplicities.
We show that this number is independent of the position of the
points and the value of the $j$-invariant and that it coincides with
the number of complex elliptic curves (with $j$-invariant $j \notin
\{0, 1728\}$).

   The result can be used to simplify G.\ Mikhalkin's algorithm to
   count curves via lattice paths (see \cite{Mi03}) in the case of rational plane
   curves.
\end {abstract}

\maketitle

\section{Introduction}
In classical algebraic geometry, the isomorphism class of an
elliptic curve is given by its $j$-invariant. The enumeration of
complex plane elliptic curves with fixed $j$-invariant using
Gromov-Witten theory has been undertaken by R. Pandharipande
\cite{Pan97}. He computed the number of elliptic curves of degree
$d$ through $3d-1$ points and with fixed $j$-invariant to be
$E(d,j)=\binom{d-1}{2}\cdot N(d)$ for $j  \notin \{0,1728\}$, where
$N(d)$ denotes the number of irreducible rational curves of
degree $d$ interpolating $3d-1$ points in general position (in case
$j\in\{0,1728\}$, the numbers differ by a factor, due to the
presence of extra automorphisms).

In tropical geometry, the isomorphism class of an elliptic curve is determined by the
(integer) length of its only cycle. This length, called tropical
$j$-invariant, can be viewed as a tropical analogue of the
$j$-invariant of an elliptic curve in classical algebraic geometry
\cite{Mi06}. The tropical reformulation of the
enumerative problem above is the enumeration of tropical curves of
genus $1$ and degree $d$ with prescribed tropical
$j$-invariant passing through a
collection of $3d-1$ points in $\R^2$ in tropical general position.

In this paper, we construct the moduli space of tropical elliptic curves as a weighted polyhedral complex. Furthermore, we define the multiplicity of a tropical elliptic curve with fixed $j$-invariant --- similar to the case of tropical rational curves --- as
the absolute value of the determinant of the evaluation map times a weight of the moduli space.
Then we define the number $E_{\trop}(d)$ to be the number of tropical elliptic curves
 with fixed tropical $j$-invariant and interpolating $3d-1$ given points in $\RR^2$, counted with multiplicities.

We prove that the numbers $E_{\trop}(d)$ are independent of the choice of the $j$-invariant for all values $j \in \R_{\geq 0}$, without exceptional values for the $j$-invariant (which is different for the complex case).
Therefore we can choose a special $j$-invariant to compute $E_{\trop}(d)$.
There are two possibilities to choose a special $j$-invariant such that we can relate the elliptic tropical curves with that $j$-invariant to rational tropical curves. One possibility is to choose a very large $j$-invariant, the other a very small $j$-invariant. We prove that an elliptic tropical curve with a very large $j$-invariant contains a contracted bounded edge in a way that its image in $\R^2$ can also be interpreted as rational tropical curve.
In this way we can show that the numbers $E_{\trop}(d)$ satisfy the equation
\begin{equation}\label{q1}E_{\trop}(d)=\binom{d-1}{2}\cdot N_{\trop}(d),\end{equation}
where $N_{\trop}(d)$ denotes the number of plane rational tropical
curves of degree $d$ through $3d-1$ points, counted with
multiplicities (see \cite{Mi03}, definition 4.15). Since by
G.\ Mikhalkin's Correspondence Theorem (see theorem 1 of \cite{Mi03})
the number of rational tropical curves $N_{\trop}(d)$ coincides with
its complex counterpart $N(d)$, it follows (using Pandharipande's
result) that $E_{\trop}(d)=E(d,j)$ for $j\not\in \{0,1728\}$. Hence
our result leads to a ``correspondence theorem'' for elliptic curves
with fixed $j$-invariant. It would be interesting to investigate
whether there is also a direct correspondence as in G.\ Mikhalkin's
theorem, that is, a bijection between the set of tropical elliptic
curves with fixed $j$-invariant (with multiplicity) and the set of
complex curves with fixed $j$-invariant (and which complex
$j$-invariant corresponds to which tropical $j$-invariant). Also, it
would be interesting to see why there is no such bijection in the
cases where the complex $j$-invariant is $j\in \{0,1728\}$.

The methods of our computation of $E_{\trop}(d)$ using a very large $j$-invariant are analogous to Pandharipande's computation of the numbers $E(d,j)$ --- we use moduli spaces of tropical elliptic curves and evaluation maps.
But we can also compute $E_{\trop}(d)$ as mentioned above in another way, using a very small $j$-invariant. Tropical curves with a very small $j$-invariant can be related to rational curves, too. Thus we can determine the number $E_{\trop}(d)$ with the aid of G.\ Mikhalkin's lattice path count (see theorem 2 of \cite{Mi03}). The computation of $E_{\trop}(d)$ using the very small $j$-invariant does not have a counterpart in complex algebraic geometry.

We think that our computation of $E_{\trop}(d)$ gives new insights in tropical geometry. As the most important example, we want to mention here the construction of the moduli space of tropical elliptic curves. This space contains cells which are equipped with weights. For some cells these weights are not natural numbers but contain a factor of $\frac{1}{2}$. This happens due to the presence of ``automorphisms''; we therefore think that our moduli space might be an example of a ``tropical orbifold''.\\
Furthermore, equating our two formulas to determine $E_{\trop}(d)$ --- the one using a very large $j$-invariant and the one using a very small $j$-invariant --- we get a new formula to enumerate tropical rational curves. Combined with G.\ Mikhalkin's lattice path algorithm to count tropical curves (see theorem 2 of \cite{Mi03}), this leads to a new lattice path count for tropical rational curves, which has the advantage that fewer paths have to be taken into account (see corollary \ref{cor-formrat}).

Note that, with some minor changes, many of our concepts can be carried over to curves
on toric surfaces other than the projective plane, as G.\ Mikhalkin's correspondence theorem
holds for these surfaces as well. But not all of our results can be generalized to these
surfaces, as e.g. corollary \ref{cor-formrat} has no counterpart for arbitrary toric surfaces, because
the notion of column-wise Newton subdivision is lost (see remarik 3.10 of \cite{GM052}).

The paper is organized as follows: in section 2, we recall some basic definitions concerning abstract and plane tropical curves and their moduli. After that, we construct in section 3 the moduli space of tropical elliptic curves as a (fractional) weighted polyhedral complex.
In section 4, we define tropical evaluation maps and use them to define multiplicities for elliptic curves with fixed $j$-invariant. Using these multiplicities, we prove in section 5 the independence of the numbers $E_{\trop}(d)$ from the configuration of the given general points and the given value of the $j$-invariant. This independence is used in section 6 to compute the numbers $E_{\trop}(d)$ and in section 7 to obtain our modified version of G.\ Mikhalkin's algorithm to count rational plane curves via lattice paths.

We would like to thank Andreas Gathmann for his inspiring ideas and for numerous helpful discussions.
\section{Tropical elliptic curves and their duals}

We will define tropical curves almost in the same way as in \cite{GM053}, with the only difference that we allow curves of higher genus. 
Let us first introduce some notions concerning graphs and metric graphs that we will need. For more details on graphs, see definition 2.1 of \cite{GM053}.
A graph can have bounded as well as unbounded edges.
We denote the set of vertices by $\Gamma^0$ and the set of edges by $\Gamma^1$. The subset of $\Gamma^1$ of bounded edges is denoted by $\Gamma^1_0$, and the subset of unbounded edges by $\Gamma^1_\infty$. Unbounded edges will also be called \emph{ends}.
A \emph{flag} $F$ of $\Gamma$ is a pair $(V,e)$ of a vertex $V$ and an edge $e$ starting at $V$. We will denote the edge $e$ of a flag $F=(V,e)$ by $[F]=e$, and the vertex $V$ by $\partial F=V$. We can think of a flag $(V,e)$ as a ``directed edge'' pointing away from the vertex $V$. For a bounded edge $e$ there are two flags $F$ and $F'$ with $[F]=e$, for an unbounded edge $e$ there is only one flag $F$ with $[F]=e$.
We denote the set of flags of $\Gamma$ by $\Gamma'$. 
Now assume $\Gamma$ is a metric graph, i.e.\ all bounded edges $e$ are equipped with a length $l(e)$ (they can be thought of as real intervals of length $l(e)$).
Given a flag $F=(V,e)$ of a bounded (respectively unbounded) edge, we can parametrize $e$ (using an affine map of slope $\pm 1$, i.e.\ a map of the form $t\mapsto c\pm t$) by an interval $[0,l(e)]$ (respectively, $[0,\infty)$), such that the vertex $V$ is at $0$. This parametrization will be called the \emph{canonical parametrization} for $F$.
The \emph{genus} of a graph $\Gamma$ is defined to be its first Betti number $g(\Gamma):= 1-\#\Gamma^0+\#\Gamma^1_0 = 1-\dim H_0(\Gamma,\Z)+\dim H_1(\Gamma,\Z) $. 

\begin{definition}
An \emph{abstract tropical curve} of genus $g$ is a metric graph $\Gamma$ of genus $g$ whose vertices have valence at least $3$. 
An \emph{abstract $n$-marked tropical curve} of genus $g$ is a tuple $(\Gamma,x_1,\ldots,x_n)$ where $\Gamma$ is an abstract tropical curve of genus $g$ and $x_1,\ldots,x_n\in \Gamma^1_\infty$ are distinct unbounded edges of $\Gamma$. We will refer to the $x_i$ as marked ends or marked points (a reason why we call them marked points is given in remark 2.7 of \cite{GM053}).
Two abstract $n$-marked tropical curves $ (\Gamma,
  x_1,\dots,x_n) $ and $ (\tilde \Gamma, \tilde x_1,\dots,\tilde x_n) $ are
  called isomorphic (and will from now on be identified) if there is a
  homeomorphism $ \Gamma \to \tilde \Gamma $ mapping $ x_i $ to $ \tilde x_i $
  for all $i$ and preserving the lengths of all bounded edges (i.e.\ every edge of $ \Gamma $ is mapped bijectively onto
  an edge of $ \tilde \Gamma $ by an affine map of slope $ \pm 1 $).

The set of all isomorphism classes of connected $n$-marked tropical curves with exactly $n$ unbounded edges and of genus $g$ is called $\mathcal{M}_{\trop,\;g,n}$.
\end{definition}

\begin{example}\label{ex-trop10}
We want to determine the space $\mathcal{M}_{\trop,\;1,1}$. An element of $\mathcal{M}_{\trop,\;1,1}$ is an abstract tropical curve with one unbounded edge, and of genus $1$. As no divalent vertices are allowed, such an abstract tropical curve consists of one bounded edge whose two endpoints are identified and glued to the unbounded edge.
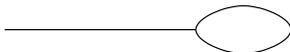
\begin{center}
\input{Graphics/absg1.pstex_t}
\end{center}
These curves only differ in the length of their bounded edge, which has to be positive. Therefore $\mathcal{M}_{\trop,\;1,1}$ is isomorphic to the open interval $(0,\infty)$.
We define $\overline{\mathcal{M}}_{\trop,\;1,1}$ to be the interval $[0,\infty)$.
Following G.\ Mikhalkin, we call the length of the bounded edge --- which is an inner
invariant of the tropical elliptic curve --- its tropical $j$-invariant, as it
plays the role of the $j$-invariant of elliptic curves in algebraic geometry (see example 3.15 of \cite{Mi06}, see also remark \ref{rem-cycle} and definition \ref{def-evj}).

\end{example}

\begin {definition} \label {def-planecurve}  
An \emph {$n$-marked plane tropical curve} of genus $g$ is
    a tuple $ (\Gamma,h,x_1,\dots,x_n) $, where $ \Gamma $ is an abstract
    tropical curve of genus $g$, $ x_1,\dots,x_n \in \Gamma^1_\infty $ are distinct
    unbounded edges of $ \Gamma $, and $ h: \Gamma \to \RR^2 $ is a continuous
    map satisfying:
    \begin {enumerate}
    \item [(i)] On each edge of $ \Gamma $ the map $h$ is of the form $ h(t) =
      a + t \cdot v $ for some $ a \in \RR^2 $ and $ v \in \ZZ^2 $ (i.e.\ ``$h$
      is affine linear with rational slope''). The integral vector $v$
      occurring in this equation if we pick for $e$ the canonical
      parametrization with respect to a chosen flag $F$ of $e$ will be denoted $ v(F) $ and called the \emph {direction} of
      $F$.
    \item [(ii)] For every vertex $V$ of $ \Gamma $ we have the \emph {balancing
      condition}
        \[ \sum_{F \in \Gamma': \partial F = V} v(F) = 0. \]
    \item [(iii)] Each of the unbounded edges $ x_1,\dots,x_n \in
      \Gamma^1_\infty $ is mapped to a point in $ \RR^2 $ by $h$ (i.e.\ $
      v(F)=0 $ for the corresponding flags).
    \end {enumerate}
    Two $n$-marked plane tropical curves $ (\Gamma,x_1,\dots,x_n,h) $ and
    $ (\tilde \Gamma, \tilde x_1,\dots,\tilde x_n,\tilde h) $ are called
    isomorphic (and will from now on be identified) if there is an isomorphism
    $ \varphi: (\Gamma,x_1,\dots,x_n) \to (\tilde \Gamma,\tilde x_1,\dots,
    \tilde x_n) $ of the underlying abstract curves such that $ \tilde h \circ \varphi = h $.
      
The \emph {degree} of an $n$-marked plane tropical curve is defined to be the
    multiset $ \Delta = \{ v(F);\; [F] \in \Gamma^1_\infty \backslash \{
    x_1,\dots,x_n \} \} $ of directions of its non-marked unbounded edges. If
    this degree consists of the vectors $ (-1,0) $, $ (0,-1) $, $ (1,1) $ each
    $d$ times then we simply say that the degree of the curve is $d$. 
 \end {definition}
\begin{remark} \label{rem-weightdir}
Note that the direction vector of a flag $F=(V,e)$ (if it is nonzero) can uniquely be written as a product of a positive integer (called the weight $\omega(e)$ of the edge $e$) and a primitive integer vector called the primitive direction $u(F)$ of the flag $F$.  
\end{remark}
\begin{definition}
For all $n\geq 0$ and $d>0$, we define $\mathcal{M}_{\trop,\;1,n}(d)$ to be the set of all isomorphism classes of connected plane tropical curves $(\Gamma,h,x_1,\ldots,x_n)$ of degree $d$ and genus $g\leq 1$.
\end{definition}

\begin{remark}\label{rem-cycle}
Note that for a connected graph $\Gamma$ the genus satisfies \begin{displaymath}g = 1-\dim H_0(\Gamma,\Z)+\dim H_1(\Gamma,\Z) = \dim H_1(\Gamma,\Z).\end{displaymath} Hence for a connected plane tropical curve of genus $1$, we have $\dim H_1(\Gamma,\Z)=1$, so there is one cycle which generates $H_1(\Gamma,\Z)$. When we avoid to pass an edge twice in both directions, there is a unique way to choose a chain of flags around this cycle (up to direction - we can go two ways around the cycle; and up to starting point). In the following, we will therefore speak of ``the cycle'' of an element of $\mathcal{M}_{\trop,\;1,n}(d)$ of genus $1$, meaning this chain of flags.
\end{remark}

The \emph{combinatorial type} of an abstract $n$-marked tropical curve $(\Gamma,x_1,\ldots,x_n)$ is
the homeomorphism class of $\Gamma$ relative $x_1,\ldots,x_n$ (i.e.\ modulo homeomorphisms that map $x_i$ to itself). We can think about it as the graph with the information of the length of the bounded edges dropped.
The \emph{combinatorial type} of a plane tropical curve $(\Gamma,h,x_1,\ldots,x_n)$ is the combinatorial type of $(\Gamma,x_1,\ldots,x_n)$ together with the directions $v(F)$ for all flags $F\in \Gamma'$.
$\mathcal{M}_{\trop,\;1,n}^{\alpha}(d)$ is defined to be the subset of $\mathcal{M}_{\trop,\;1,n}(d)$ of tropical curves of combinatorial type $\alpha$.
\begin{definition}
Let $\alpha$ be a combinatorial type in the space $\mathcal{M}_{\trop,\;1,n}(d)$. The \emph{deficiency} $\defi(\alpha)$ is defined to be
 \begin{equation*}
  \defi(\alpha)=
    \begin{cases}
   2 \mbox{ if $g=1$ and the cycle is mapped to a point in $\mathbb{R}^2$},\\
    1 \mbox{ if $g=1$ and the cycle is mapped to a line in $\mathbb{R}^2$},\\
    0 \mbox{ otherwise}.
    \end{cases}
 \end{equation*}
We will also speak of the deficiency $\defi(C)$ of a curve $C$.
\end{definition}
Having defined elliptic tropical curves we now want to come to their dual Newton subdivisions.

Let $V$ be an $r$-valent vertex of a plane tropical curve
$(\Gamma,h)$ (without markings) and let $e_1,\ldots, e_r$ be the
counterclockwise enumerated edges adjacent to $V$. Draw in the
$\Z^2$-lattice an orthogonal line $L(e_i)$ of integer length
$\omega(e_i)$ (where $\omega(e)$ denotes the weight of $e$, see remark \ref{rem-weightdir}) to
$h(e_i)$, where $L(e_1)$ starts at any lattice point and $L(e_i)$
starts at the endpoint of $L(e_{i-1})$, and where by ``integer
length'' we mean $\# (\Z^2 \cap L(e_i)) -1$. The balancing condition
tells us that we end up with a closed $r$-gon. If we do this for
every vertex we end up with a polygon in $\Z^2$ that is divided into
smaller polygons. The polygon is called the \emph{Newton polygon} of
the tropical curve, and the division the corresponding \emph{Newton
subdivision}. Note that the ends of the curve correspond to
line segments on the boundary of the Newton polygon. The
Newton polygon of a curve of degree $d$ is the triangle $\Delta_d$
with vertices $(0,0)$, $(0,d)$ and $(d,0)$. For more details on the
dual Newton subdivision of a tropical curve, see \cite{Mi03},
section 3.4.

Some properties of plane tropical curves can be read off from their dual picture. Here are some examples:
  \begin{enumerate}
\item A plane tropical curve is called \emph{simple} if its dual subdivision contains only triangles
    and parallelograms (see definition 4.2 of \cite{Mi03}). (This property can also be defined without using the dual language.)
  \item The \emph{genus} of a simple plane tropical curve $ (\Gamma,h) $ is
    equal to the number of lattice points of the
    subdivision contained in the interior of the Newton polytope minus the number of parallelograms (see lemma 4.6 of \cite{Mi03}).
  \item \label{tropdef-e}
    Let $V$ be a trivalent vertex of $\Gamma$ and $e_1,e_2,e_3$ the edges
    adjacent to $V$. The \emph{multiplicity} of $V$ is defined to be the area of the parallelogram spanned by the two directions of $e_1$ and $e_2$. (Due to the balancing condition this is independent of the choice of $e_1$ and $e_2$.) It is equal to $2$ times the
    area of the dual triangle (see definition 2.16 of \cite{Mi03}).
  \item \label{tropdef-f}
    The \emph{multiplicity} $\mult(C)$ of a $3$-valent tropical plane curve is the product over all multiplicities of the vertices. In the dual language, the multiplicity of a simple curve is the product over all double areas of triangles of the dual subdivision (see definition 4.15 of \cite{Mi03}).
  \end{enumerate}
For more details on dual Newton subdivisions, see for example \cite{Ma06}, section 5.
\section{The moduli space of tropical elliptic curves}
Let us study the space $\mathcal{M}_{\trop,\;1,n}(d)$.
There are only finitely many combinatorial types in the space $\mathcal{M}_{\trop,\;1,n}(d)$ (analogously to 2.10 of \cite{GM053}).
\begin{lemma}\label{lem-dimofmalpha}
For every combinatorial type $ \alpha $ occurring in $\mathcal{M}_{\trop,\;1,n}(d)$ the space $\mathcal{M}^{\alpha}_{\trop,\;1,n}(d)$ is naturally an (unbounded) open convex polyhedron
  in a real vector space of dimension $2+ \#\Gamma^1_0$, that is a subset of a real vector space given by finitely many linear equations and
  finitely many linear strict inequalities.
The dimension of $\mathcal{M}^{\alpha}_{\trop,\;1,n}(d)$ is equal to
\begin{displaymath}
\dim (\mathcal{M}^{\alpha}_{\trop,\;1,n}(d))=3d+n+g-1-\sum_{V \in \Gamma^0}(\val V-3)+\defi(\alpha)
\end{displaymath}
\end{lemma}
\begin{proof}
Fixing a combinatorial type means we fix the homeomorphism class of $\Gamma$ and the directions of all flags. We do not fix the lengths $l(e)$ of the bounded edges. (Note that the length of an image $h(e)\subset \R^2$ is determined by $l(e)$ and the direction.) Also, we can move an image $h(\Gamma)$ in the whole plane.
Choose a root vertex $V$ of $\Gamma$. Two coordinates are given by the position of the image $h(V)$ in the plane. $\#\Gamma^1_0$ coordinates are given by the lengths of the bounded edges (which have to be positive).
Hence we can embed $\mathcal{M}_{\trop,\;1,n}^{\alpha}(d)$ in $\R^{2+\#\Gamma^1_0}$.
Note that a graph with $3d+n$ unbounded edges and of genus $1$ has $3d+n-\sum_{V\in \Gamma^0}(\val V-3)$ bounded edges, whereas a rational graph with $3d+n$ unbounded edges has $3d+n-3-\sum_{V\in \Gamma^0}(\val V-3)$ bounded edges (see for example \cite{Mi03}, proof of 2.13).
For a genus $1$ curve $C$, the lengths of the bounded edges are not independent however, as some are contained in the cycle. So these lengths satisfy two conditions, namely that their images have to close up a cycle in $\R^2$. These conditions are only independent if $\defi(C)=0$. If $\defi(C)=1$ there is one independent equation, and if $\defi(C)=2$ there is none. Hence the dimension of the polyhedron $\mathcal{M}_{\trop,\;1,n}^{\alpha}(d)$ in $\R^{2+\#\Gamma^1_0}$ is $2+\#\Gamma^1_0-2+\defi(C)= 3d+n-\sum_{V\in \Gamma^0}(\val V-3)+\defi(C)$. Note that the polyhedron is unbounded, because we can for example move the image of the root vertex in the whole plane.
For rational curves, no equations have to be fulfilled and we have $ \dim \mathcal{M}_{\trop,\;1,n}^{\alpha}(d)= 3d+n-1-\sum_{V\in \Gamma^0}(\val V-3)$.
\end{proof}
\begin{proposition}\label{rem-boundarymalpha}
Let $\alpha$ be a combinatorial type occurring in $\mathcal{M}_{\trop,\;1,n}(d)$. Then every point in $\overline{\mathcal{M}}^{\alpha}_{\trop,\;1,n}(d)$ (where the closure is taken in $\R^{2+\#\Gamma^1_0}$, see lemma \ref{lem-dimofmalpha}) can naturally be thought of as an element in  $\mathcal{M}_{\trop,\;1,n}(d)$. The corresponding map \begin{displaymath} i_\alpha: \overline{\mathcal{M}}^{\alpha}_{\trop,\;1,n}(d) \rightarrow \mathcal{M}_{\trop,\;1,n}(d) \end{displaymath} maps the boundary $
  \partial \mathcal{M}^{\alpha}_{\trop,\;1,n}(d) $ to a union of strata $\mathcal{M}^{\alpha'}_{\trop,\;1,n}(d) $ such that
  $\alpha'$ is a combinatorial type with fewer internal edges than $\alpha$.
  Moreover, the restriction of $ i_\alpha $ to any inverse image of such a
  stratum $ \mathcal{M}^{\alpha'}_{\trop,\;1,n}(d)  $ is an affine map.
\end{proposition}
\begin{proof}
Note that by the proof of \ref{lem-dimofmalpha} a point in the boundary of the open polyhedron $\mathcal{M}^{\alpha}_{\trop,\;1,n}(d)\subset \R^{2+\#\Gamma^1_0}$ corresponds to a tuple $(\Gamma,h,x_1,\ldots,x_n)$ where some edges $e$ have length $l(e)=0$.
Such a curve is of a different combinatorial type then, because the homeomorphism class of the graph has changed. For all edges $e$ with length $l(e)=0$ the vertices $\partial F$ and $\partial F'$ of the two flags $F$ and $F'$ with $[F]=[F']=e$ are identified. We can as well remove the edges of length $0$ then.
Note that the balancing condition will be fulfilled at the new vertices.
Two examples what this can look like are shown in the following picture. The edges which tend to have length zero when we move towards the boundary of the open polyhedron $\mathcal{M}^{\alpha}_{\trop,\;1,n}(d)$ are drawn in bold.
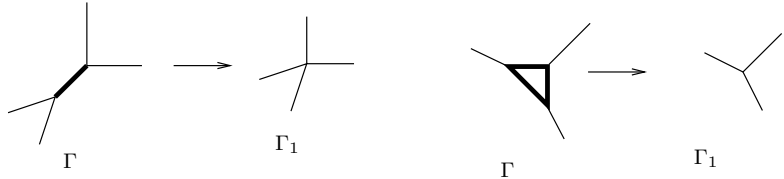
\begin{center}
\input{Graphics/boundary3.pstex_t}
\end{center}
Let $\Gamma_1$ be the graph which is obtained by removing the edges of length $0$. Note that $\Gamma_1$ has fewer bounded edges than $\Gamma$. The tuple $(\Gamma_1,h|_{\Gamma_1},x_1,\ldots,x_n)$ is a tropical curve again, possibly of a smaller genus than $(\Gamma,h,x_1,\ldots,x_n)$. This shows that the points in
  the boundary $ \partial\mathcal{M}^{\alpha}_{\trop,\;1,n}(d) $ can naturally be thought of
  as parametrized tropical curves in $\mathcal{M}_{\trop,\;1,n}(d) $ themselves. The combinatorial types $\alpha'$ that can occur in the boundary of $\mathcal{M}^{\alpha}_{\trop,\;1,n}(d)$, that is, in the image $ i_\alpha(\partial \mathcal{M}^{\alpha}_{\trop,\;1,n}(d)) $, have by construction fewer bounded edges than $\alpha$. Finally, it is clear that the restriction of $ i_\alpha $ to the inverse
  image of any stratum $\mathcal{M}^{\alpha'}_{\trop,\;1,n}(d)  $ is an affine map since
  the affine structure on any stratum is given by the position of the curve in
  the plane and the lengths of the bounded edges.
\end{proof}
\begin{definition}\label{def-containedinbound}
We will say that a type $\alpha'$ appears \emph{in the boundary} of another type $\alpha$, if there is a point in $\partial \mathcal{M}^{\alpha}_{\trop,\;1,n}(d)$ that is identified with a curve of type $\alpha'$ (as in the proof of proposition \ref{rem-boundarymalpha}).
\end{definition}
Our aim is to define a slightly different moduli space
$\mathcal{M}'_{\trop,\;1,n}(d)$, where the strata of dimension
bigger than $3d+n$ are excluded and where we add a \emph{weight} for
each stratum. In fact, our moduli space should be something similar
to an abstract tropical variety, a weighted polyhedral complex:
\begin{definition}
Let $ X_1,\dots,X_N $ be (possibly unbounded) open convex polyhedra in real
  vector spaces. A \emph {polyhedral complex} with cells $ X_1,\dots,X_N $ is a
  topological space $X$ together with continuous inclusion maps $ i_k:
  \overline {X_k} \to X $ such that $X$ is the disjoint union of the sets $
  i_k(X_k) $ and the ``coordinate changing maps'' $ i_k^{-1} \circ i_l $ are
  linear (where defined) for all $ k \neq l $. We will usually drop the
  inclusion maps $ i_k $ in the notation and say that the cells $ X_k $ are
  contained in $X$.

  The \emph {dimension} $ \dim X $ of a polyhedral complex $X$ is the maximum
  of the dimensions of its cells. We say that $X$ is of \emph {pure dimension}
  $ \dim X $ if every cell is contained in the closure of a cell of dimension
  $ \dim X $. A point of $X$ is said to be \emph {in general position} if it is
  contained in a cell of dimension $ \dim X $. For a point $P$ in general position, we denote the cell of dimension $\dim X$ in which it is contained by $X_P$.

A \emph{weighted polyhedral complex} is a polyhedral complex such that there is a weight $\omega(X_i)\in \Q$ associated to each cell $X_i$ of highest dimension.
\end{definition}
We are now ready to define the moduli space $\mathcal{M}'_{\trop,\;1,n}(d)$, which is important for our methods:
\begin{definition}\label{def-relevant}
Remove the strata of dimension bigger than $3d+n$ from $\mathcal{M}_{\trop,\;1,n}(d)$. Also, remove the strata of rational curves which are not contained in the boundary of a genus $1$ curve as in \ref{def-containedinbound}.
Let $\alpha$ be a type such that $\dim \mathcal{M}^{\alpha}_{\trop,\;1,n}(d)=3d+n$.
We associate the following weights to the strata of dimension $3d+n$:
\begin{enumerate}
\item Assume $\defi(\alpha)=0$, and the curves of type $\alpha$ are of genus $1$.
As we have already seen, the condition that the image of the cycle closes up in $\R^2$ is given by two independent linear equations $a_1$ and $a_2$ on the lengths of the bounded edges. We associate as weight the index of the map
\begin{displaymath}\binom{a_1}{a_2}: \Z^{2+\#\Gamma^1_0}\rightarrow \Z^2.
\end{displaymath}
(For more details on lattices, maps between vector spaces and lattices and their indices, see \cite{Rau06}.)\label{def-rel-1}
\item Assume $\defi(\alpha)=1$. Assume first that the $4$-valent vertex is adjacent to the cycle, that is, locally the curves look like the following picture:
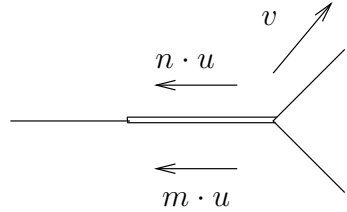
\begin{center}
\input{Graphics/flatloop.pstex_t}
\end{center}
In the notations above, $n\cdot u$, $m\cdot u$ and $v$ denote the direction vectors of the corresponding edges ($n$ and $m$ are chosen such that their greatest common divisor is $1$).
If $n\neq m$, or if $n=m=1$ and the cycle is formed by three edges due to the presence of a marked point, we associate the weight $|\det(u,v)|$. \label{def-rel-2}
If $n=m=1$ and no point is on the flat cycle, then we associate $\frac{1}{2}|\det(u,v)|$. (Due to the balancing condition this definition is not dependent of the choice of $v$.)\\
In case the $4$-valent vertex is not adjacent to the cycle, we associate the weight $0$.
\item Assume $\defi(\alpha)=2$. Assume first that the $5$-valent vertex is adjacent to the cycle, that is, locally the curves look like this:
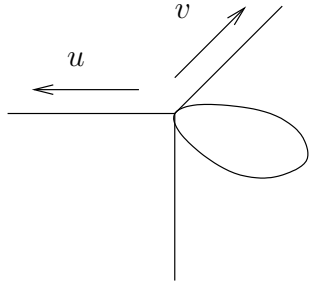
\begin{center}
\input{Graphics/five.pstex_t}
\end{center}
where $u$ and $v$ denote the direction vectors of the corresponding edges. We associate the weight $\frac{1}{2}(|\det(u,v)|-1)$. (Note that due to the balancing condition this definition is independent of the choice of $u$ and $v$.)
In the case that there are two $4$-valent vertices or that the $5$-valent vertex is not adjacent to the cycle, we associate the weight $0$.
\label{def-rel-3}
\end{enumerate}
The strata of dimension $3d+n$ or less together with these weights for the strata of top dimension form the space $\mathcal{M}'_{\trop,\;1,n}(d)$, called the \emph{moduli space of (relevant) elliptic tropical curves}.
\end{definition}
 The reason to drop the cells of dimension bigger than $3d+n$ is that we want to construct later on a morphism to a polyhedral complex of the same dimension $3d+n$. The strata of dimension bigger than $3d+n$ would not be mapped injectively to the image. We will only be interested in strata which are mapped injectively, therefore the strata of dimension higher than $3d+n$ are not important to us and we can drop them.  
\begin{remark}
Note that the definitions of weight do not depend on the choice of coordinates for the cell $\mathcal{M}^{\alpha}_{\trop,\;1,n}(d)$. This is clear for each case except the first one. In the first case, the two equations given by the cycle do not depend on the choice of a root vertex, they depend only on the choice of an order for the bounded edges. But this corresponds to an isomorphism on $\R^{2+\#\Gamma^1_0}$ of determinant $\pm 1$, therefore the index of $\binom{a_1}{a_2}$ does no depend on this choice.

The definition of the weights for the different strata seems somewhat unnatural. However, we will see that the weights are just the right ones for our proofs later on. The idea behind this definition is that we think of elliptic tropical curves as rational tropical curves with two additional marked points, whose images we require to coincide.
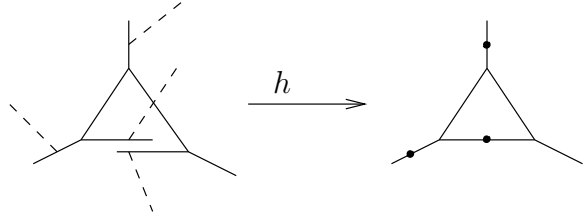
\begin{center}
\input{Graphics/elr.pstex_t}
\end{center}
The space of rational curves with two additional marked points is of course bigger than the space of elliptic curves, but it contains the space of elliptic curves as the kernel of the map $a_1\times a_2$.
Note that the index of $a_1\times a_2$ is equal to
\begin{displaymath}\gcd(a_1)\cdot \gcd(a_2)\cdot \chi\big(\ker(a_1),\ker(a_2)\big),\end{displaymath} where $\chi\big(\ker(a_1),\ker(a_2)\big)$ denotes the index of the sublattice generated by $\ker(a_1)+\ker(a_2)$ in $\Z^{2+\#\Gamma1_0}$ (see example 1.5 of \cite{Rau06}).

If $\defi(C)=1$, the weight we choose is derived in the same way: we compute the multiplicity of a rational curve with two additional marked points.
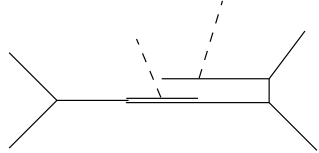
\begin{center}
\input{Graphics/elr1.pstex_t}
\end{center}
Sometimes we included a factor of $\frac{1}{2}$ in our weights. It seems maybe unnatural to allow weights which are not natural numbers. But the factors of $\frac{1}{2}$ are only necessary when the cycle of the elliptic curve ``allows automorphisms'': when it is a loop consisting of one edge with two non distinguishable orientations, or when it consists of two non distinguishable edges. Due to these factors we believe that the moduli space we construct here can be thought of as a ``tropical orbifold''.
\end{remark}
\begin{lemma}\label{lem-tildempoly}
 The space $\mathcal{M}'_{\trop,\;1,n}(d)$ (defined in \ref{def-relevant}) is a weighted polyhedral complex of pure dimension $3d+n$.
\end{lemma}
\begin{proof}
The cells are obviously the strata $\mathcal{M}^{\alpha}_{\trop,\;1,n}(d)$ corresponding to relevant types. By \ref{lem-dimofmalpha} they are open convex polyhedra. By proposition \ref{rem-boundarymalpha}, their boundary is also contained in $\mathcal{M}'_{\trop,\;1,n}(d)$, and the coordinate changing maps are linear.
By definition, the highest dimension of a relevant cell is $3d+n$. Furthermore, by definition each rational type which is contained in $\mathcal{M}^{\alpha}_{\trop,\;1,n}(d)$ is in the boundary of a type of genus $1$. Each higher-valent vertex can be resolved to $3$-valent vertices. Therefore each type is contained in the boundary of a type of codimension $0$. By definition, the strata of top dimension are equipped with weights as required.
\end{proof}
\section{The multiplicity of an elliptic tropical curve}\label{sec-mult}
We also want to define morphisms between weighted polyhedral complexes. 

\begin {definition} \label{def-morph}
    A \emph {morphism} between two weighted polyhedral complexes $X$ and $Y$ is a
    continuous map $ f:X \to Y $ such that for each cell $ X_i \subset X $ the
    image $ f(X_i) $ is contained in only one cell of $Y$, and $ f|_{X_i} $ is
    a linear map (of polyhedra). 

Assume $f: X \to Y $ is a morphism of weighted polyhedral complexes of the same pure
    dimension, and $ P \in X $ is a point such that both $P$ and $ f(P) $
    are in general position (in $X$ resp.\ $Y$).
Then locally around $P$ the
    map $f$ is a linear map between vector spaces of the same dimension. We
    define the \emph {multiplicity} $ \mult_f(P) $ of $f$ at $P$ to be the
    absolute value of the determinant of this linear map times the weight of the cell $X_P$. Note that the
    multiplicity depends only on the cell $X_P$ of $X$ in which $P$ lies. We will
    therefore also call it the multiplicity of $f$ in this cell.

 A point $ Q \in Y $ is said to be \emph {in $f$-general
    position} if $Q$ is in general position in $Y$ and all points of $
    f^{-1}(Q) $ are in general position in $X$. Note that the set of points in
    $f$-general position in $Y$ is the complement of a subset of $Y$ of
    dimension at most $ \dim Y-1 $; in particular it is a dense open subset.
    Now if $ Q \in Y $ is a point in $f$-general position we define the \emph
    {degree} of $f$ at $Q$ to be
      \[ \deg_f(Q) := \sum_{P \in f^{-1}(Q)} \mult_f(P). \]
    Note that this sum is indeed finite: first of all there are only finitely
    many cells in $X$. Moreover, in each cell (of maximal dimension) of $X$
    where $f$ is not injective (i.e.\ where there might be infinitely many
    inverse image points of $Q$) the determinant of $f$ is zero and hence so is
    the multiplicity for all points in this cell.

    Moreover, since $X$ and $Y$ are of the same pure dimension, the cones of
    $X$ on which $f$ is not injective are mapped to a locus of codimension at
    least 1 in $Y$. Thus the set of points in $f$-general position away from
    this locus is also a dense open subset of $Y$, and for all points in this
    locus we have that not only the sum above but indeed the fiber of $Q$ is
    finite.
 \end {definition}

Note that the definition of multiplicity $\mult_f(P)$ in general depends on the coordinates we choose for the cells. However, we will use this definition only for the morphism $\ev \times j$ (see \ref{def-evj}) for which the absolute value of the determinant does not depend on the chosen coordinates, if they are chosen in a natural way (in our case this means we choose a lattice basis of the space $\mathcal{M}^{\alpha}_{\trop,\;1,n}(d)\subset \R^{1+\#\Gamma^1_0}$, see remark \ref{rem-detunique}).

The following maps will be important to count elliptic curves:
\begin{definition}\label{def-evj}
Let
\begin{displaymath}\ev_i:\mathcal{M}'_{\trop,\;1,n}(d)\rightarrow \R^2,\; (\Gamma,h,x_1,\ldots ,x_n) \mapsto h(x_i)\end{displaymath} denote the $i$-th evaluation map. By $\ev=\ev_1\times \ldots\times \ev_n$ we denote the combination of all $n$ evaluation maps.

If $C$ is an elliptic curve, let $ \Gamma_1 $ be the minimal connected subgraph of genus $1$ of $ \Gamma $ that
    contains the unbounded edge $ x_1 $. Note that  $\Gamma_1 $
    cannot contain vertices of valence 1. So if we ``straighten'' the graph $
    \Gamma_1 $ at all 2-valent vertices (that is we replace the two adjacent edges
    and the vertex by one edge whose length is the sum of the lengths of the
    original edges) then we obtain an element of $\overline{\mathcal{M}}_{\trop,\;1,1}=[0,\infty) $ that we denote by
    $ j(C) $. (For an elliptic curve, $j(C)\neq 0$.)
If $C$ is rational, we define $j(C)=0\in \overline{\mathcal{M}}_{\trop,\;1,1}=[0,\infty) $.
We call $j(C)$ the \emph{tropical $j$-invariant} of $C$.

A combination of these maps yields
\begin{displaymath}\ev \times j: \mathcal{M}'_{\trop,\;1,n}(d)\rightarrow \R^{2n}\times \overline{\mathcal{M}}_{\trop,\;1,1}.
\end{displaymath}
\end{definition}
\begin{example}
The following picture shows an elliptic curve $C$. The marked points are drawn as dotted lines. The subgraph $\Gamma_1$ is indicated with a bold dotted line. The image $j(C)$ is an abstract tropical curve where the cycle has length $l_1+l_2+\ldots+l_8$.
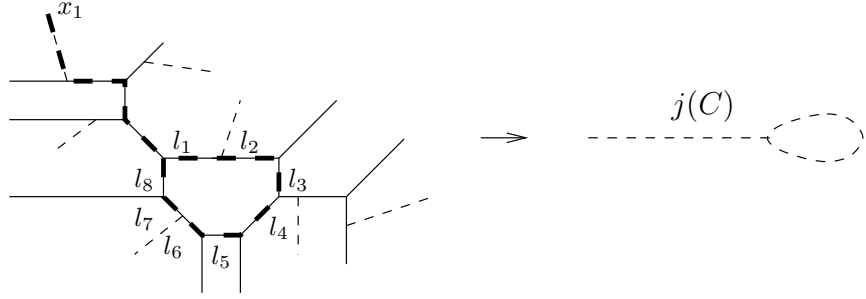
\begin{center}
\input{Graphics/jex.pstex_t}
\end{center}
\end{example}
\begin{lemma}\label{lem-evlinear}
The map $ev\times j$ restricted to a stratum $\mathcal{M}^{\alpha}_{\trop,\;1,n}(d)$ is a linear map.
\end{lemma}
\begin{proof}
The coordinates on $\mathcal{M}^{\alpha}_{\trop,\;1,n}(d)$ are by \ref{lem-dimofmalpha} given by a root vertex $V$ and an order on the bounded edges. Of course, these coordinates do not need to be independent, but if they are not, they fulfill a linear condition themselves.
As $\Gamma$ is connected, we can reach $x_i$ from the root vertex $V$ by a chain of flags $F$, such that $[F]$ is a bounded edge. Then the position of $h(x_i)$ is given as a sum \begin{displaymath}h(V)+\sum_F v(F)\cdot l([F]),\end{displaymath} where the summation goes over all flags $F$ in the chain. Hence the position $h(x_i)$ is given by two linear expressions in the coordinates of $\mathcal{M}^{\alpha}_{\trop,\;1,n}(d)$.
The length of the bounded edge of $j(C)$ is by definition given as the sum of the lengths of all bounded edges contained in the cycle of $\Gamma_1$ that we straightened to get $j(C)$.
\end{proof}
\begin{lemma}\label{lem-evpolymap}
Let $n=3d-1$. Then the map $\ev\times j $ is a morphism of weighted polyhedral complexes of the same dimension.
\end{lemma}
\begin{proof}
The space $\mathcal{M}'_{\trop,\;1,n}(d)$ is a weighted polyhedral complex of dimension $3d+n=3d+3d-1=6d-1$ by \ref{lem-tildempoly}. The space $\R^{2n}\times  \overline{\mathcal{M}}_{\trop,\;1,1}$ is by \ref{ex-trop10} isomorphic to $\R^{2n}\times [0,\infty)$. This is obviously a polyhedral complex with only one cell, and we can make it weighted by associating the weight $1$ to this cell. Its dimension is $2n+1=2(3d-1)+1=6d-1$, too.
The map $\ev\times j$ restricted to a cell of the weighted polyhedral complex $\mathcal{M}'_{\trop,\;1,n}(d)$ is linear by \ref{lem-evlinear}. Furthermore, it maps each cell into the one cell of the space $\R^{2n}\times  \overline{\mathcal{M}}_{\trop,\;1,1}$.
\end{proof}
We will from now on assume that $n=3d-1$, in order to have a morphism between polyhedral complexes of the same dimension.
\begin{remark}
Fix $n$ points $p_1,\ldots,p_n$ and a $j$-invariant $l$ in $\ev\times j$-general position (see definition \ref{def-morph}). Then determine the set of tropical elliptic curves $(\ev\times j)^{-1}(p_1,\ldots,p_n,l)$ which pass through the points and have $j$-invariant $l$. We count each such elliptic curve with its $\ev\times j$-multiplicity, that is, we determine $\deg_{\ev\times j}((p_1,\ldots,p_n,l))$.
Define $E_{\trop}(d)=\deg_{\ev\times j}((p_1,\ldots,p_n,l))$ to be the \emph{number of tropical elliptic curves through $3d-1$ points in general position and with fixed $j$-invariant}. Our aim is to show that this definition does not depend on the choice of $(p_1,\ldots,p_n,l)$. This statement will be shown in theorem \ref{thm-degconst}.
\end{remark}
\begin{remark}\label{rem-detunique}
Given a stratum $\mathcal{M}^{\alpha}_{\trop,\;1,n}(d)$ of top dimension we choose a lattice basis of the space
$\mathcal{M}^{\alpha}_{\trop,\;1,n}(d)$ in $\R^{2+\#\Gamma^1_0}$, which contains the lattice $\Z^{2+\#\Gamma^1_0}$.
(For $\R^{2+\#\Gamma^1_0}$, we choose natural coordinates given by a root vertex and an order of the bounded edges.)
With this choice, we can compute a matrix representation of $\ev \times j$ and hence compute its determinant.
We claim that the absolute value of this determinant (and thus, the $\ev \times j$-multiplicity) does not depend on the choices we made.
To see this, note first that a different choice of the root vertex or the order of the bounded edges corresponds to a basis change of determinant $\pm 1$ (see remark 3.2 of \cite{GM053}).
If we choose a different lattice basis for $\mathcal{M}^{\alpha}_{\trop,\;1,n}(d)$, then we have to multiply the matrix of $\ev \times j$ with the basis change matrix. But as this is a basis change of lattice bases, it is of determinant $\pm 1$ and does therefore not change the absolute value of the determinant of $\ev \times j$.
\end{remark}
The following remark helps us to determine $\mult_{C}(\ev\times j)$ in some cases:
\begin{remark}\label{rem-johannes}
Let $\alpha$ be a type of $3$-valent genus $1$ curves with
$\defi(\alpha)=0$. We want to determine the multiplicity of $\ev
\times j$ in the stratum $\mathcal{M}^{\alpha}_{\trop,\;1,n}(d)$. By
definition, it is equal to the absolute value of the determinant of $\ev \times j$ times
the weight of the cell $\mathcal{M}^{\alpha}_{\trop,\;1,n}(d)$. The
weight of the cell is defined to be the index of the lattice map
$\binom{a_1}{a_2}$, where $a_1$ and $a_2$ denote the two equations
given by the cycle (see \ref{def-relevant}). To compute the $\ev
\times j$-multiplicity, we need to compute this weight, and then a
matrix representation of $\ev \times j$ restricted to $\ker(a_1)\cap
\ker(a_2)$. To get this matrix representation, we need a lattice
basis of the subspace $\mathcal{M}^{\alpha}_{\trop,\;1,n}(d)\subset
\R^{2+\#\Gamma^1_0}$. However, lattice bases are in general not easy
to determine. We can instead use example 1.7 of \cite{Rau06}. It
states that $|\det(\ev\times j)|$ times the lattice index is equal
to the absolute value of the determinant of the map
\begin{displaymath}\ev \times j \times a_1 \times
a_2:\R^{2+\#\Gamma^1_0}\rightarrow \R^{6d-2}\times
\overline{\mathcal{M}}_{\trop,\;1,1} \times \R^2\end{displaymath}
(where $\ev$ and $j$ denote here matrix representations of $\ev$
respectively $j$ in the coordinates given by the root vertex and the
lengths of all bounded edges). This map goes from the space
$\R^{2+\#\Gamma^1_0}$ which surrounds
$\mathcal{M}^{\alpha}_{\trop,\;1,n}(d)$. More precisely,
$\mathcal{M}^{\alpha}_{\trop,\;1,n}(d)$ is equal to $\ker(a_1)\cap
\ker(a_2)\subset \R^{2+\#\Gamma^1_0}$. The bigger space
$\R^{2+\#\Gamma^1_0}$ does not parametrize tropical curves, as the
length coordinates of a general vector of $\R^{2+\#\Gamma^1_0}$ do
not need to fulfill the conditions given by the cycle. But (after
choosing coordinates and chains of flags from the root vertex to
each marked point) we can still write down the matrix of $\ev \times
j \times a_1 \times a_2$. Note that while the matrix does depend on
the choices we make, the absolute value of the determinant does not
(see remark 4.65 of \cite{Ma06}). That is, to compute the
multiplicity of $\ev \times j$ in
$\mathcal{M}^{\alpha}_{\trop,\;1,n}(d)$, we can write down a matrix
representation for the map $\ev \times j \times a_1 \times a_2$ and
compute its determinant. In the following, we will denote the map by
$\ev \times j \times a_1 \times a_2$, even though it is not uniquely
determined by this term. It depends on the chosen matrix
representation, that is, on the chosen chains of flags. We will keep
in mind that $|\det(\ev \times j \times a_1 \times a_2)|$ is
uniquely determined, though.
\end{remark}
\begin{example}
Compute $|\det(\ev \times j \times a_1 \times a_2)|$ for the following (local) picture of an elliptic curve.
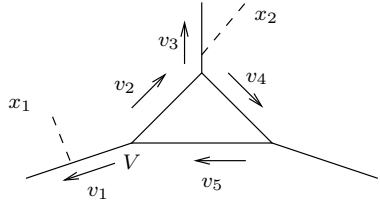
\begin{center}
\input{Graphics/evj.pstex_t}
\end{center}
Let $V$ be the root vertex.
We choose the following chains of flags: for $x_1$, we pass $v_1$. For $x_2$, we pass $v_2$ and $v_3$. For $a_1\times a_2$, we pass $v_2$, $v_4$ and $v_5$. Then the matrix reads:

\[ \left( \begin {array}{cccccc}
         E_2&v_1&0&0&0&0 \\
         E_2&0&v_2&v_3&0&0 \\
         0&0&v_2&0&v_4&v_5\\
0&0&1&0&1&1
       \end {array} \right)\]

Each row except the last represent two rows, the first column represents two columns. $E_2$ stands for the $2\times 2$ unit matrix.
\end{example}
We add four other statements that help to determine the $\ev\times j$-multiplicity in some cases:
\begin{lemma}\label{lem-evjcontr-1}
Let $\defi(C)=2$, that is, $C$ contains a contracted loop. $C'$ denotes the rational curve which arises from $C$ if we remove the loop.

If the vertex $V$ to which the loop is adjacent is $5$-valent, then \begin{displaymath}\mult_{\ev\times j}(C)= \frac{1}{2}\big(\mult(V)-1\big) \cdot \mult(C'),\end{displaymath} where $\mult(V)$ denotes the multiplicity of the vertex $V$ of $C'$ from which the loop was removed and $\mult(C')$ denotes the multiplicity of the tropical curve $C'$.

Else the $\ev\times j$-multiplicity of $C$ is $0$.
\end{lemma}
\begin{proof}
To determine a matrix representation for $\ev\times j$, we do not have to consider equations given by the cycle and lattice bases --- the lengths of the bounded edges are independent. That is, we can choose a root vertex and an order of the bounded edges and write down the matrix of $\ev\times j$ with respect to this basis. Note that in the $j$-row there is just one unit at the coordinate of the contracted edge, as no other edge is contained in the cycle. To compute the determinant, we can therefore remove this last line and the column of the coordinate of the contracted edge which forms the cycle. The remaining matrix consists of the evaluation maps in the $3d-1$ marked points, and it does not take the contracted edge into account. That is, this matrix describes the evaluation in the marked points of the rational curve $C'$ which arises when we remove the contracted edge. Proposition 3.8 of \cite{GM053} tells us that the determinant of this matrix is equal to the multiplicity of the rational curve $C'$. To determine the $\ev\times j$-multiplicity, we also have to multiply with the weight of the stratum in which $C$ lies.
The only case where the weight is non-zero is \ref{def-relevant}\ref{def-rel-3}, so we can assume that the contracted loop is adjacent to a $5$-valent vertex $V$, and denote the direction vectors of two of the other adjacent edges by $u$ and $v$.
Then the weight of the stratum is $\frac{1}{2}(|\det(u,v)|-1)$ by definition.
If we remove the loop and consider the rational curve $C'$, then this weight is equal to $\frac{1}{2}(\mult(V)-1)$, where $\mult(V)$ denotes the multiplicity of the vertex $V$ to which the contracted loop was adjacent.
\end{proof}
\begin{lemma}\label{lem-evjcontr-2}
Let $C$ be a ($3$-valent) curve with a contracted bounded edge $e$, which is not a loop. Let $C'$ denote the rational curve that arises if we remove $e$ and straighten the two $2$-valent vertices $V_1$ and $V_2$ emerging like this.
Let $u$ denote the direction of a remaining edge adjacent to $V_1$, and let $v$ denote the direction of a remaining edge adjacent to $V_2$.

If $e$ is part of the cycle, then \begin{displaymath}\mult_{\ev\times j}(C)= |\det(u,v)|\cdot \mult(C').\end{displaymath}

Else the $\ev\times j$-multiplicity of $C$ is $0$.
\end{lemma}
\begin{proof}
We assumed that the contracted bounded edge $e$ is adjacent to two different $3$-valent vertices. (These vertices have to be $3$-valent, as we only compute the $\ev\times j$-multiplicity in strata of top dimension.)

The balancing condition implies that at each of these two vertices the two other adjacent edges are mapped to opposite directions:
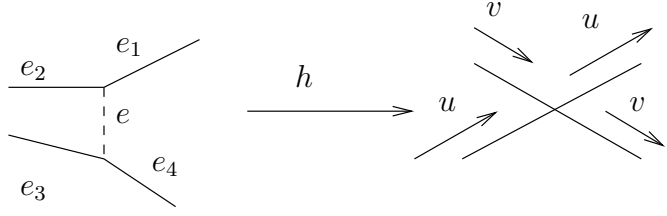
\begin{center}
\input{Graphics/contr.pstex_t}
\end{center}
We are going to use remark \ref{rem-johannes} to compute the $\ev\times j$-multiplicity, that is, we use a matrix representation of $\ev\times j \times a_1\times a_2$, where $a_1$ and $a_2$ denote the two equations of the cycle.
If $e$ is not part of the cycle its length is a coordinate which is needed neither for the evaluations nor for the map $j$, so we get determinant $0$. Assume now that $e$ is part of the cycle.
Then $e$ is only needed for the map $j$, as it is contracted by $h$ it is not needed to describe a position of the image of a marked point. That is, in the column for $e$, we have a $1$ at the row of the map $j$ and $0$ in every other row.
In the matrix of $\ev\times j \times a_1\times a_2$, we can choose a chain of flags to each marked point which avoids the edge $e$. This is possible, as $e$ is contained in the cycle. In the two equations $a_1$ and $a_2$, exactly one of the edges $e_1$ and $e_2$ will take part, and also exactly one of the edges $e_3$ and $e_4$. Assume without loss of generality that $e_1$ and $e_3$ are part of the cycle. If $e_1$ is used in a chain of flags to a marked point, then also $e_2$, and if $e_3$, then also $e_4$. Let $l_i$ denote the length of $e_i$, and $l$ denote the length of $e$.
Assume that the directions are as labelled in the picture.
Then the $\ev\times j \times a_1\times a_2$-matrix looks like this:

\[ \begin {array}{l|ccccccc}
         & h(V) & l_1 & l_2 & l_3 & l_4 & l&\mbox{other edges} \\  \hline
         \mbox {marked points using neither of the $e_i$} &
           E_2 & 0 & 0 & 0 & 0 & 0 &* \\
         \mbox {marked point using $e_1$} &
           E_2 & u & u & 0 & 0 & 0 &* \\
         \mbox {marked points using $e_3$} &
           E_2 & 0 & 0 & v & v& 0 &*\\
         a_1,a_2 &
           0 & u & 0 & v & 0 & 0 &* \\
         \mbox {coordinate of $\overline{\mathcal{M}}_{\trop,\;1,1}$} &
           0 & 1 & 0 & 1 & 0 & 1&*
       \end {array} \]

We perform the following operations which do not change the absolute value of the determinant:
we delete the last row and the $l$-column. We subtract the $l_2$-column from the $l_1$-column and we change the place of the $l_1$-column: it shall appear as first column. Then, we subtract the $l_4$-column from the $l_3$-column and move the $l_3$-column to the second place. At last, we put the two rows $a_1$ and $a_2$ to the beginning.
After these operations the matrix looks like this:

 \[ \begin {array}{l|cc|cccc}
        &    l_1 & l_3 & h(V) & l_2 & l_4 &\mbox{other edges} \\  \hline
a_1,a_2 &
         u & v & 0 & 0 & 0 &*  \\ \hline
\mbox {marked points using neither of the $e_i$} &
           0 & 0 & E_2 & 0 & 0 &* \\
\mbox {marked point using $e_1$} &
           0 & 0 & E_2& u&0&* \\
 \mbox {marked points using $e_3$} &
           0 & 0 & E_2 & 0 & v &*
       \end {array} \]

Note that this matrix is now a block matrix with a $2\times 2$ block on the top left, and a block that we will denote by $A$ on the bottom right. That is, its determinant is equal to $\det(u,v)\cdot \det A$. Now $A$ is the matrix of the evaluation map in the $3d-1$ marked points of the rational curve $C'$ which arises from $C$ when we remove the contracted bounded edge $e$ and straighten the two bounded edges $e_1$ and $e_2$ as well as the two bounded edges $e_3$ and $e_4$ to one edge.
As before, proposition 3.8 of \cite{GM053} tells us that $|\det A|=\mult(C')$. That is, the $\ev\times j$-multiplicity of $C$ is equal to $|\det(u,v)|$ times the multiplicity of the rational curve $C'$ which arises after removing the contracted edge and straightening the adjacent edges.
Our argument here assumes that the edges $e_1,\ldots,e_4$ are all bounded. However, we can prove the same if some of these edges are not bounded. Their lengths do not appear as coordinates then, but also there cannot be marked points behind unbounded edges.
\end{proof}
\begin{lemma}\label{lem-evjcontr-3}
Let $C$ be a curve with a flat cycle (that is, $\defi(C)=1$) which is not adjacent to a marked point. If there is no $4$-valent vertex adjacent to the flat cycle, then $\mult_{\ev\times j}(C)=0$. Else, with the notations as in the picture below, we have
\begin{displaymath}\mult_{\ev \times j}(C)=\begin{cases} (m+n)\cdot |\det(u,v)|\cdot \mult(C')\mbox{ if } m \neq n \\ |\det(u,v)|\cdot \mult(C') \mbox{ if } m=n=1\end{cases}\end{displaymath}
where $C'$ denotes the rational curve that arises if we glue the two edges that form the cycle to one edge of direction $(m+n)\cdot u$ and straighten the $2$-valent vertex emerging like this.
\end{lemma}
\begin{proof}
The following picture shows the flat cycle of the curve $C$. We choose $m$ and $n$ such that $\gcd(m,n)=1$.
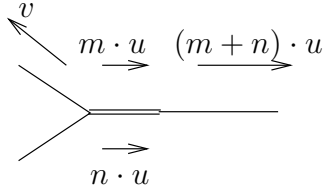
\begin{center}
\input{Graphics/fl.pstex_t}
\end{center}
To determine the matrix of $\ev\times j$, we need a lattice basis of $\mathcal{M}^{\alpha}_{\trop,\;1,n}(d)$. As the equations of the cycle are given by $l_1\cdot m\cdot u-l_2\cdot n\cdot u$, we can choose unit vectors for all coordinates except $l_1$ and $l_2$, plus the vector with $n$ at the $l_1$-coordinate and $m$ at the $l_2$-coordinate.
As $\gcd(m,n)=1$, this is a lattice basis. The $j$-invariant of $C$ is given by $l_1+l_2$. That is, in the $j$-row of the matrix, we have only zeros except for the column which belongs to the vector with $n$ at $l_1$ and $m$ at $l_2$, there we have the entry $m+n$. But then we can delete the $j$-row and this column. The determinant we want to compute is equal to $(m+n)$ times the determinant of the matrix which arises after deleting. This matrix can easily be seen to be the matrix of evaluating the points of the rational curve $C'$ which arises after identifying the two edges which form the cycle. Due to \cite{GM053} proposition 3.8, its determinant is equal to $\mult(C')$. The factor of $|\det(u,v)|$ (respectively, $\frac{1}{2}\cdot |\det(u,v)|$ if $n=m=1$) has to be included, because this is by definition \ref{def-relevant} \ref{def-rel-2} the weight of the stratum $\mathcal{M}^{\alpha}_{\trop,\;1,n}(d)$.
\end{proof}
\begin{lemma}\label{lem-evjcontr-4}
Let $C$ be a curve with $\defi(C)=0$ and such that the cycle is formed by $3$ edges.

Then $\mult_{\ev \times j}(C)= \mult(V)\cdot \mult(C')$, where $C'$ is the rational curve which arises when we shrink the cycle to a vertex $V$.
\end{lemma}
\begin{proof}
We compute the $\ev\times j$-multiplicity of $C$ using remark \ref{rem-johannes}. That is, we compute the determinant of the matrix $\ev\times j \times a_1\times a_2$, where $a_1$ and $a_2$ denote the two equations of the cycle.
The following picture shows the curve locally around the cycle and fixes a labeling of the adjacent edges:
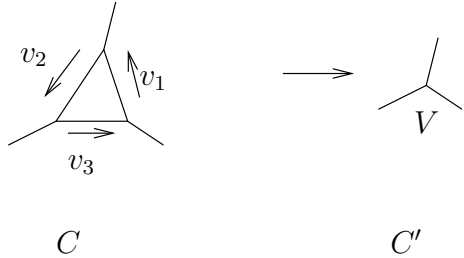
\begin{center}
\input{Graphics/dreieck.pstex_t}
\end{center}
Note that the matrix $\ev\times j \times a_1\times a_2$ has a block form with a $0$ block on the bottom left, because the equations of the cycle and the $j$-invariant only need the three length coordinates of $v_1$, $v_2$ and $v_3$.
The block on the top left is just the evaluation of the rational curve $C'$ at the marked points --- hence by \cite{GM053} proposition 3.8, its determinant is equal to $\mult(C')$. So $\mult_{\ev \times j}(C)$ is equal to $\mult(C')$ times the absolute value of the determinant of the matrix

\[ \left( \begin {array}{ccc}
        1&1&1\\
        v_1&v_2&v_3
       \end {array}\right) \]

where the last row stands for the two rows given by the equation of the cycle. The absolute value of this determinant can be computed to be $|\det(v_1,v_2)|+|\det(v_1,v_3)|+|\det(v_2,v_3)|$ which is --- using the dual picture, for example --- easily seen to be equal to $\mult V$:
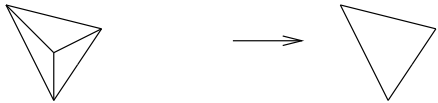
\begin{center}
\input{Graphics/dickesdreieck.pstex_t}
\end{center}
The sum of the three determinants is equal to the double areas of the three small triangles, $\mult V$ is equal to the double area of the big triangle.
\end{proof}
\section{The number of tropical elliptic curves with fixed $j$-invariant}
A \emph{string} in a tropical curve $C$ is a subgraph of $\Gamma$ homeomorphic either to $\R$ or to $S^1$ (that is, a ``path'' starting and ending with an unbounded edge, or a path around a loop) that does not intersect the closures $\overline{x_i}$ of the marked points (see also \cite{Ma06}, definition 4.47).
If the number of marked points on $C$ is less than $3d+g-1$, then $C$ has a string. This follows from lemma 3.50 of \cite{Ma06}. We will need the notion of a string in the following proof of the main theorem of this section:
\begin{theorem}\label{thm-degconst}
Let $n=3d-1$. The degrees $\deg_{\ev \times j}(\mathcal{P})$ do not depend on $\mathcal{P}$. (Here $\mathcal{P}=(p_1,\ldots,p_{3d-1},l)\in \R^{2n}\times \overline{\mathcal{M}}_{\trop,\;1,1}$ denotes a configuration in $\ev\times j$-general position consisting of $3d-1$ points in $\R^2$ and a length $l$ for the $j$-invariant.)
\end{theorem}
\begin{proof}
Analogously to the proof of 4.4 of \cite{GM053}, we have that the degree of $ \ev \times j$ is \emph {locally} constant on the
  subset of $ \R^{6d-2} \times \overline{\mathcal{M}}_{\trop,\;1,1} $ of points in $ \ev \times j $-general position
  since at any curve that counts for $ \deg_{\ev\times j}(\mathcal{P}) $ with a non-zero
  multiplicity the map $ \ev \times j $ is a local isomorphism. The points
  in $ \ev\times j $-general position are the complement of a polyhedral complex of
  codimension 1, that is they form a finite number of top-dimensional regions
  separated by ``walls'' that are polyhedra of codimension 1. Hence it remains to show that $ \deg_{\ev\times j} $ is locally constant at
  these points, too. Such a general point on a wall is the image under $
  \ev \times j $ of a general tropical curve $C$ of a type $\alpha$ such that $\mathcal{M}^{\alpha}_{\trop,\;1,n}(d)$ is of codimension 1. So we have to check that $ \deg_{\ev \times j} $ is locally
  constant around such a point $ C \in \mathcal{M}'_{\trop,\;1,n}(d)$.
More precisely, if $\mathcal{P}$ is such a point on a wall, and $C$ is a curve through $\mathcal{P}$, we want to show that the sum of the $\ev\times j$-multiplicities of the curves through $\mathcal{P}'$ near $\mathcal{P}$ and close to $C$ does not depend on $\mathcal{P}'$.
Let us determine what types $\alpha$ are of codimension $1$, using \ref{lem-dimofmalpha}.
\begin{enumerate}
\item \label{case1}$\defi(\alpha)=0$, $\alpha$ is of genus $1$ and has one $4$-valent vertex (besides the $3$-valent vertices);
\item\label{case2} $\defi(\alpha)=1$ and $\alpha$ has two $4$-valent vertices;
\item \label{case3}$\defi(\alpha)=1$ and $\alpha$ has one $5$-valent vertex;
\item \label{case4}$\defi(\alpha)=2$ and $\alpha$ has three $4$-valent vertices;
\item \label{case5}$\defi(\alpha)=2$ and $\alpha$ has one $5$-valent and one $4$-valent vertex;
\item \label{case6}$\defi(\alpha)=2$ and $\alpha$ has one $6$-valent vertex.
\end{enumerate}
Note that the codimension $1$ case that $\alpha$ is the type of a rational curve is missing here: the reason is that we do not ``cross'' such a wall consisting of rational curves, we can only enlarge the $j$-invariant if $j=0$, not make it smaller. More precisely, the curves which pass through a configuration $\mathcal{P}'$ in the neighborhood of a point configuration through which a rational curve passes, are always of the same types; the types (and with them, the multiplicities with which we count) do not depend on $\mathcal{P}'$.

For each of the cases in the list, we have to prove separately that $ \deg_{\ev \times j} $ is locally
  constant around a curve $C$ of type $\alpha$.
The proof for \ref{case1} is similar to the proof of 4.4 in \cite{GM053}.
There are three types which have $\alpha$ in their boundary. The following is a local picture:
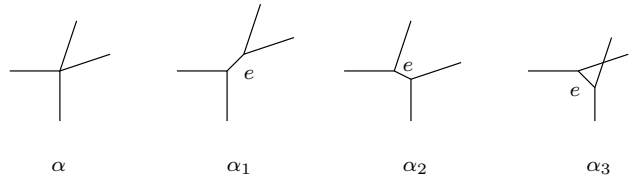
\begin{center}
\input{Graphics/four.pstex_t}
\end{center}
To compute the $\ev\times j$-multiplicity of a curve of each type, we use remark \ref{rem-johannes}, that is, we compute a matrix representation of $\ev\times j\times a_1\times a_2$, where $a_1,a_2$ denote the two equations of the cycle. We can choose the coordinates in such a way that these three matrices only differ in one column - in the column corresponding to the new edge $e$. Then we can use the same operations as in 4.4 of \cite{GM053} to prove that the sum of the three determinants is $0$. The matrices we use here differ from the ones in 4.4 of \cite{GM053} because they contain the two lines corresponding to $a_1$ and $a_2$, and the line corresponding to the $j$-invariant. However, the argument does not change in the presence of these other lines.
Also, with the same argument as in 4.4 of \cite{GM053}, the question whether there is a curve of type $\alpha_i$ through a given configuration $\mathcal{P}'$ close to $\mathcal{P}$ depends on the sign of the determinant. So we can conclude that we either get the types where the determinant has positive sign or the types where it has negative sign. But as the sum is $0$, the sum of the absolute values of those determinants, for which a curve exists, stays constant and does not depend on $\mathcal{P}'$.

So let us now come to \ref{case2}.
First note that if none of the $4$-valent vertices is adjacent to the flat cycle, then we count all curves of a type which has $\alpha$ in its boundary with the weight $0$, so we do not have to consider this case.
So at least one of the $4$-valent vertices is adjacent to the flat cycle.
If exactly one of the $4$-valent vertices is adjacent to the flat cycle, then the only curves which have $\alpha$ in their boundary and which do not count with weight $0$ are the curves where the other $4$-valent vertex is resolved, as in case \ref{case1}. The proof is then analogous to the proof of case \ref{case1}, only using the matrices of $\ev \times j$ instead of the big matrices of $\ev \times j \times a_1\times a_2$.
So we can assume now that both $4$-valent vertices are adjacent to the flat cycle.
Assume first that none of the edges adjacent to a $4$-valent vertex in the flat cycle is a marked point.
We claim that $C$ has a string.
Consider the connected components of $\Gamma \backslash \bigcup_i \overline {x_i}$.
As in the proof of 3.50 of \cite{Ma06}, remove the closures of the marked points $\overline {x_1},\ldots,\overline {x_n}$ from $\Gamma$ one after the other. We only remove edges at $3$-valent vertices. Therefore each removal can either separate one more component, or break a cycle. Assume that all connected components are rational (else $C$ contains a string). Then one of our removals must have broken the cycle. As $C$ is marked by $3d-1$ points, we end up with $3d-1$ connected components. But then there has to be one connected component which contains two unbounded edges, hence $C$ contains a string.

 If $C$ has at least two strings then $C$
 moves in an at least 2-dimensional family with the images of the marked points fixed.  As $ \mathcal{M}_{\trop,\;1,1} $ is one-dimensional this
  means that $C$ moves in an at least 1-dimensional family with the image point
  under $j $ fixed. But then also the curves close to $C$ are not fixed, hence they count $0$. So we do not have to consider this case. Also, if for all curves $C'$ which contain $C$ in their boundary the string does not involve an edge of the cycle, then $C$ (and all curves $C'$) move in an at least 1-dimensional family with the image point
  under $j $ fixed. So we do not have to consider this case either.

So we assume now that $C$ lies in the boundary of a type which has exactly one string that involves (at least) one of the edges of the flat cycle.

There are (up to symmetry) five possibilities for the string. We will show them in the following local picture.

Assume now that there is a marked point adjacent to a $4$-valent vertex of the flat cycle. Then the removal of this marked point both breaks a cycle and separates two components. So we cannot conclude that $C$ has a string. However, we can conclude that there is no other marked point adjacent to the cycle, as else two marked points would map to the same line. Hence in this case the curve looks locally like our sixth picture below.
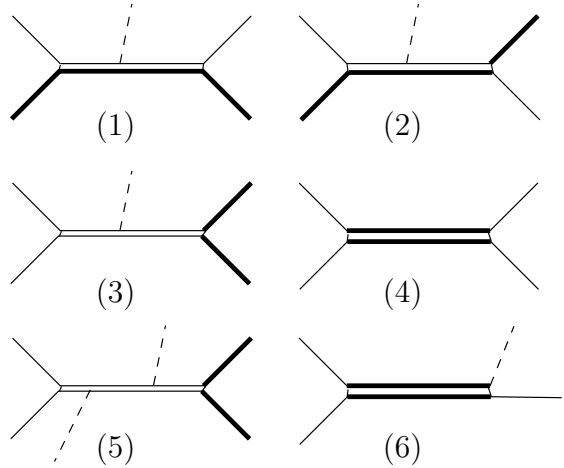
\begin{center}
\input{Graphics/string14.pstex_t}
\end{center}

As the cycle is not a string in the cases (1)-(3) and (5), there must be a marked point adjacent to it. Two marked points adjacent to the flat cycle are only possible if the string does not involve any edge of the flat cycle (as in (3) and (5)).

In each of the six cases, there are four types which contain $\alpha$ in their boundary (see \ref{rem-boundarymalpha}). The following picture shows the four types $\alpha_1,\ldots,\alpha_4$ for case (1).
We will give our argument only for case (1), it is analogous in all other five cases.

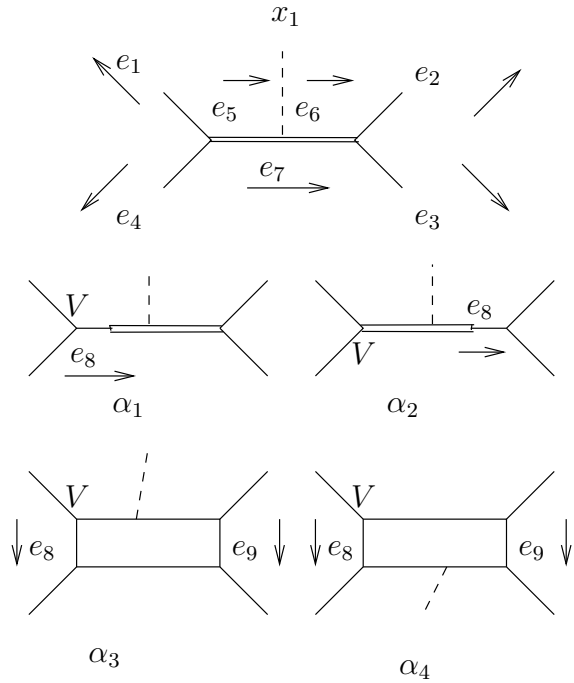
\begin{center}
\input{Graphics/4types1.pstex_t}
\end{center}

\begin{center}
\input{Graphics/4types.pstex_t}
\end{center}

Our first aim is to show that we can choose bases of the corresponding strata $\mathcal{M}^{\alpha}_{\trop,\;1,n}(d)$ such that the matrices $\ev\times j$ (and the matrices $\ev\times j \times a_1\times a_2$ for the types $\alpha_3$ and $\alpha_4$) contain a block which involves only the edges locally around the flat cycle. We can then make statements about the $\ev\times j$-multiplicity (with the aid of \ref{rem-johannes}) using these smaller blocks.

Choose the root vertex for the four types to be $V$ as indicated in the picture above. Also choose the labeling for the edges around the cycle as above.
Let $v_i$ be the directions of the $e_i$ (as indicated in the picture). We then have $v_5=v_6$. As $e_5$, $e_6$ and $e_7$ are mapped to the same line in $\R^2$, we can choose $n$ and $m$ with $v_5=v_6=n\cdot u$ and $v_7=m\cdot u$ such that $\gcd(n,m)=1$.
We will consider a matrix representation $A_i$ ($i=1,2$) of $\ev \times j$ for the types $\alpha_1$ and $\alpha_2$ and a matrix representation $B_i$ ($i=3,4$) of $\ev \times j \times a_1\times a_2$ for $\alpha_3$ and $\alpha_4$.
The $\ev \times j$-multiplicity for $\alpha_1$ is then given by $|\det(u,v_2)\cdot \det A_1|$, the $\ev \times j$-multiplicity for $\alpha_2$ is $|\det(u,v_1)\cdot \det A_2|$ due to \ref{def-relevant}\ref{def-rel-2}. For $\alpha_3$, it is due to remark \ref{rem-johannes} given by $|\det B_3|$ and for $\alpha_4$ by $|\det B_4|$. Later on, we will also need to consider matrix representations $A_3$ and $A_4$ of $\ev \times j$ for the types $\alpha_3$ and $\alpha_4$. We will however not choose a lattice bases for those, so they are not useful for the computation of the $\ev \times j$-multiplicity. We will specify later on what bases we choose for $A_3$ and $A_4$.

We choose a basis of the subspace $\mathcal{M}^{\alpha_i}_{\trop,\;1,n}(d)\subset \R^{2+\#{\Gamma_i}_{1}^0}$ for $i=1,2$ consisting of two unit vectors for the root vectors and unit vectors for all bounded edges except $e_5$, $e_6$ and $e_8$. In addition, we take two vectors with $e_5$, $e_6$ and $e_7$-coordinates as follows: $(1,-1,0)$ and $(0,m,n)$. In fact, this is a lattice basis:
As $\gcd(n,m)=1$, we can find integer numbers such that $am+bn=1$. Then we can complete our basis with the vector $(0,b,a)$ (at $e_5,e_6,e_7$) and get a lattice basis of $\Z^{2+\#{\Gamma_i}_{1}^0}$.
For $i=3,4$, we choose a basis of $\mathcal{M}^{\alpha_i}_{\trop,\;1,n}(d)\subset \R^{2+\#{\Gamma_i}_{1}^0}$ consisting of only unit vectors except three vectors involving the coordinates of $e_5,\ldots,e_9$.

Because the bases we choose for the $A_i$ and for the $B_i$ differ only by a few vectors, there will be a block in which the matrices $A_i$ and $B_i$ ($i=3,4$) do not differ. The following argumentation works therefore analogously for all six matrices $A_1,\ldots,A_4$, $B_3$ and $B_4$.

Assume that $d_1$ unbounded (nonmarked) edges can be reached from $V$ via $e_1$, $d_2$ via $e_2$ and so on. As the only string passes via $e_4$ and $e_3$, there must be $d_1$ marked points which can be reached from $V$ via $e_1$, $d_2$ marked points via $e_2$, $d_3-1$ marked points via $e_3$ and $d_4-1$ via $e_4$.
Note that the marked points which can be reached via $e_1$ and $e_2$ do not need any of the length coordinates of edges via $e_3$ or $e_4$. As there are $2\cdot (d_3-1+d_4-1)$ rows for the marked points via $e_3$ and $e_4$ and $2d_3-1 +1 -3 +1+ 2d_4-2$ bounded edges via $e_3$ and $e_4$, all six matrices have a $0$ block on the top right.
For $B_3$ and $B_4$, we also put the equations $a_1$ and $a_2$ of the cycle in the first block of rows.

\[ \begin {array}{l|c|c|c}
         & h(V) & \mbox{other edges}& \mbox{ edges via }e_3 \mbox{ and }e_4   \\ \hline & & & \\

         \mbox {$x_1$, pts behind $ e_1 $ and $e_2$ and $j$-coord} &
           E_2 & * & 0  \\
         \mbox {pts behind $ e_3 $ and $e_4$} &
           E_2 & * & * \\

       \end {array} \]

The block on the bottom right is the same for all six matrices. So we can disregard it and only consider the top left block given by the marked points which can be reached via $e_1$ and $e_2$, the $j$-coordinate, $a_1$ and $a_2$ for $B_3$ and $B_4$ and the length coordinates of $e_1,e_2,e_5,\ldots e_8/e_9$ plus the length coordinates of bounded edges via $e_1$ respectively $e_2$. Choose a marked point $x_2$ which can be reached via $e_1$ and a marked point $x_3$ which can be reached via $e_2$.
Choose the following order for the rows: begin with the marked points $x_1,\ldots,x_3$, then take the $j$-coordinate (and for the matrices $B_3$ and $B_4$ the equations $a_1,a_2$ of the cycle). Then take the remaining marked points.
Choose the following order for the columns: begin with $h(V)$, $e_1$ and $e_2$. For the types $\alpha_1$ and $\alpha_2$, take the two basis vectors involving $e_5,\ldots,e_7$ and then $e_8$. For the types $\alpha_3$ and $\alpha_4$, take the three basis vectors involving $e_5,\ldots,e_9$. For $B_3$ and $B_4$, take the length coordinates of $e_5,\ldots,e_9$. Then take the remaining length coordinates.
Note that each marked point which can be reached via $e_1$ has the same entries in the first 7 (respectively, 9 for $B_3$ and $B_4$) columns as $x_2$. Each marked point which can be reached via $e_2$ has the same entries in the first 7 (respectively, 9) columns as $x_3$. That is, we can subtract the $x_2$-rows from all rows of marked points via $e_1$ and the $x_3$-rows from all rows of marked points via $e_2$. Then we have a $0$ block on the bottom left. Note that the bottom right block is equal for all six matrices. That is, we can now go on with the four $7\times 7$-matrices and the two $9\times 9$ matrices. The determinants of the original six matrices only differ by the factor which is equal to the determinants of the corresponding $7\times 7$-matrices (respectively, $9\times 9$) matrices.
Let us call these blocks $A_i'$, respectively $B_i'$.
Here are the four blocks $A_1'$, $A_2'$, $B_3'$ and $B_4'$ and their determinants:

\[ A_1'=\left( \begin {array}{cccccc}
         E_2&0&0&n\cdot u&0&(n+m)\cdot u \\
        E_2&v_1&0&0&0&0 \\
         E_2&0&v_2&0&n\cdot m\cdot u&(n+m)\cdot u\\
0&0&0&0&m+n&0
       \end {array} \right)\]

\begin{displaymath}
\det(A_1')=-n\cdot (n+m)^2\cdot \det(u,v_1)\cdot \det(u,v_2)
\end{displaymath}

\[ A_2'=\left( \begin {array}{cccccc}
         E_2&0&0&n\cdot u&0&0  \\
         E_2&v_1&0&0&0&0\\
         E_2&0&v_2&0&n\cdot m\cdot u&(n+m)\cdot u\\
0&0&0&0&m+n&0
       \end {array} \right)\]

\begin{displaymath}
\det(A_2')=-n\cdot (n+m)^2\cdot \det(u,v_1)\cdot \det(u,v_2)
\end{displaymath}

\[ B_3'=\left( \begin {array}{cccccccc}
         E_2&0&0&n\cdot u&0&0&0&0\\
         E_2&v_1&0&0&0&0&0&0\\
         E_2&0&v_2&n\cdot u&n\cdot  u&0&0&0\\
0&0&0&n\cdot u&n\cdot u&-m\cdot u& v_1+n\cdot u&-v_2+n\cdot u\\
0&0&0&1&1&1&1&1
       \end {array} \right)\]

\begin{displaymath}
\det(B_3')=-\det(u,v_1)\cdot \det(u,v_2)\cdot n\cdot \big( (n^2+nm)(\det(u,v_1)+\det(u,v_2))+n\det(v_1,v_2)\big)
\end{displaymath}

\[ B_4'=\left( \begin {array}{cccccccc}
          E_2&0&0&n\cdot u&0&0&-v_1-m\cdot u&0\\
          E_2&v_1&0&0&0&0&0&0\\
         E_2&0&v_2&n\cdot u&n\cdot  u&0&-v_1-m\cdot u&v_2-m\cdot u\\
0&0&0&n\cdot u&n\cdot u&-m\cdot u& -v1-m\cdot u&v_2-m\cdot u\\
0&0&0&1&1&1&1&1
       \end {array} \right)\]

\begin{displaymath}
\det(B_4')=\det(u,v_1)\cdot \det(u,v_2)\cdot n\cdot \big( (m^2+nm)(\det(u,v_1)+\det(u,v_2))-n\det(v_1,v_2)\big)
\end{displaymath}
A computation shows that
\begin{equation}\label{eq1}
\det(u,v_2)\cdot \det A_1'+ \det(u,v_1) \cdot \det A_2'- \det B_3+\det B_4=0.\end{equation}

Note that in the cases (4) and (6) \emph{without marked points adjacent to the flat cycle} we have to make a difference if $n=m=1$. In this case, definition \ref{def-relevant}\ref{def-rel-2} tells us that we have to multiply the types analogous to $\alpha_1$ and $\alpha_2$ (which still contain a flat cycle) with the factor $\frac{1}{2}\cdot \det(u,v_1)$ respectively $\frac{1}{2}\cdot \det(u,v_2)$, instead of $\det(u,v_1)$ and $\det(u,v_2)$.
However, the types $\alpha_3$ and $\alpha_4$ are not different in these cases, so we count them only once. Altogether, the weighted sum of determinants as above is still $0$.

We still need to check which types occur for a given point configuration $\mathcal{P}''$ near $\mathcal{P}'$.
Let $\mathcal{P}''\subset \R^{2n}\times \overline{\mathcal{M}}_{\trop,\;1,1}$ be a configuration. If there exists a curve $C$ of type $\alpha_i$ through $\mathcal{P}''$, then $A_i^{-1}\cdot \mathcal{P}''$ gives us the coordinates of $C$ in $\mathcal{M}^{\alpha_i}_{\trop,\;1,n}(d)$ in the basis $\{v_{i,1},\ldots,v_{i,2n-1}\}$. That it, the first two coordinates of the vector \begin{displaymath}\sum_j (A_i^{-1}\cdot \mathcal{P}'')_j \cdot v_{i,j}\subset \mathcal{M}^{\alpha_i}_{\trop,\;1,n}(d)\end{displaymath} denote the position of the root vertex, and all other coordinates the lengths of the bounded edges of $C$. A curve of type $\alpha_i$ exists if and only if all coordinates of the vector $\sum_j (A_i^{-1}\cdot \mathcal{P}'')_j \cdot v_{i,j}\subset \mathcal{M}^{\alpha_i}_{\trop,\;1,n}(d)$ which correspond to lengths are positive.
Choose $\mathcal{P}''$ close to the configuration $\mathcal{P}'$, through which a curve of type $\alpha$ exists. By continuity of $A_i^{-1}$, all coordinates of $\sum_j (A_i^{-1}\cdot \mathcal{P}'')_j \cdot v_{i,j}\subset \mathcal{M}^{\alpha_i}_{\trop,\;1,n}(d)$ except the length of $e_8$ ($i=1,2$), respectively of $e_8$ and $e_9$ ($i=3,4$), are positive.

Note that there is a curve of type $\alpha_i$ ($i=1,2$) through $\mathcal{P}''$ if and only if the $e_8$-coordinate of $A_i^{-1}\cdot \mathcal{P}''$ is positive.

Now we specify which bases we choose for the types $\alpha_3$ and $\alpha_4$. For $\alpha_3$, begin again with the two unit vectors for the position of the root vertex. Take unit vectors for all bounded edges which are not contained in the cycle.
Let \begin{align*}&M_1:=-\det(v_1,v_2)+n\cdot \det(v_1,u)-n\cdot \det(u,v_2), \\ &M_2:=-n\cdot \det(u,v_2) \mbox{ and } \\& M_3:=n\cdot \det(u,v_1).\end{align*} Take the three vectors with entries \begin{align*}&(-1,1,0,0,0),\\& (0,m,n,0,0) \mbox{ and } \\&(0,M_1,0,-M_2,M_3)\end{align*} at the coordinates of $e_5,\ldots,e_9$. (Let the vector with the entries $(0,M_1,0,-M_2,M_3)$ be the last basis vector.) These three vectors are linearly independent and satisfy the conditions given by the cycle. However, we cannot say whether this basis is a lattice basis of $\mathcal{M}^{\alpha_3}_{\trop,\;1,n}(d)$. So we do not know whether the determinant of the matrix $A_3$ is equal to the $\ev \times j$ multiplicity of $\alpha_3$. But we are not interested in the determinant of $A_3$ here, we just want to use $A_3$ to check whether there is a curve of type $\alpha_3$ through $\mathcal{P}''$ or not.
Note that the last basis vector is the only one which involves the lengths of $e_8$ and $e_9$. As due to the balancing condition we have $\det(u,v_1)>0$ and $\det(u,v_2)>0$, the two entries of this vector corresponding to these two lengths are positive.
That is, there is a curve of type $\alpha_3$ through $\mathcal{P}''$ if and only if the last coordinate of $A_3^{-1}\cdot \mathcal{P}''$ (that is, the coordinate with which we have to multiply our last basis vector to get the lengths) is positive.

For $\alpha_4$, choose besides the unit vectors the three vectors with entries \begin{align*}&(-1,1,0,0,0), \\& (0,m,n,0,0) \mbox{ and } \\& (0,M_1',0,M_2,-M_3)\end{align*} at the coordinates of $e_5,\ldots,e_9$, where \begin{displaymath}M_1':=-\det(v_1,v_2)-m\cdot \det(u,v_1)-m\cdot \det(u,v_2).\end{displaymath} The two entries of $(0,M_1',0,M_2,-M_3)$ corresponding to the lengths of $e_8$ and $e_9$ are negative.
Hence there is a curve of type $\alpha_4$ through $\mathcal{P}''$ if and only if the last coordinate of $A_4^{-1}\cdot \mathcal{P}''$ (that is, the coordinate with which we have to multiply our last basis vector to get the lengths) is negative.

For all four types, we are interested in the last coordinate of $A_i^{-1}\cdot \mathcal{P}''$. By Cramer's rule, this last coordinate is equal to $\det \tilde{ A_i} /\det A_i$, where $\tilde A_i$ denotes the matrix where the last column of $A_i$ is cancelled and replaced by $\mathcal{P}''$. Note that the four matrices $A_1,\ldots,A_4$ only differ in the last column. Hence the matrices $\tilde{A_i}$ do not depend on $i$, and we can decide whether there is a curve of type $\alpha_i$ through $\mathcal{P}''$ by determining the sign of $\det A_i$. (This argument is analogous to the proof of proposition 4.4 in \cite{GM053}.)

Recall that $|\det A_i|$ is a product of a factor which does not differ for all four types and a factor which is equal to the determinant of a $7\times 7$-matrix $A_i'$ which describes a curve of type $\alpha_i$ ``locally around the cycle''.

Here are the two matrices $A_3'$ and $A_4'$ and their determinants:

\[ A_3'=\left( \begin {array}{cccccc}
          E_2&0&0&n\cdot u&0&0 \\
          E_2&v_1&0&0&0&0\\
         E_2&0&v_2&0&n\cdot m\cdot u&-M_2\cdot (-v_1-n\cdot u)+M_3\cdot (v_2-n\cdot u)\\
0&0&0&0&m+n&M_1-M_2+M_3
       \end {array} \right)\]

\begin{displaymath}
\det(A_3')=\det(u,v_1)\cdot \det(u,v_2)\cdot n^2\cdot \big( (n^2+nm)(\det(u,v_1)+\det(u,v_2))+n\det(v_1,v_2)\big)
\end{displaymath}

\[ A_4'=\left( \begin {array}{cccccc}
         E_2&0&0&n\cdot u&0&M_2\cdot(-v_1-m\cdot u)\\
          E_2&v_1&0&0&0&0\\
         E_2&0&v_2&0&n\cdot m\cdot u&0\\
0&0&0&0&m+n&M_1'+M_2-M_3
       \end {array} \right)\]

\begin{displaymath}
\det(A_4')=-\det(u,v_1)\cdot \det(u,v_2)\cdot n^2\cdot \big( (m^2+nm)(\det(u,v_1)+\det(u,v_2))-n\det(v_1,v_2)\big)
\end{displaymath}
We know that $\det(u,v_1)\geq 0$, $\det(u,v_2)\geq 0$ and $\det(v_1,v_2)\geq0$.
So there are now two cases to distinguish:
\begin{itemize}
\item $( (m^2+nm)(\det(u,v_1)+\det(u,v_2))-n\det(v_1,v_2)\big)\geq0$ --- then $\det A_4'$ is negative. As we have seen, a curve of type $\alpha_4$ exists if and only if the last coordinate of $A_4'^{-1}\cdot \mathcal{P}''$ is negative, hence if and only if $\det \tilde{A_i}$ (a matrix which depend only on $\mathcal{P}''$, not on $i$) is positive.
$\det A_1'$, $\det A_2'$ are both negative, a curve of one of these types exists if the last coordinate of $A_i'^{-1}\cdot \mathcal{P}''$ ($i=1,2$) is positive, hence if $\det \tilde{A_i}$ is negative.
$\det A_3'$ is positive, and a curve of this type exists if the last coordinate of $A_3'^{-1}\cdot \mathcal{P}''$ is positive, hence if $\det \tilde{A_i}$ is positive.
Hence $\alpha_1$ and $\alpha_2$ are on one side of the ``wall'', $\alpha_3$ and $\alpha_4$ on the other. But as
in this case equation \ref{eq1} from above reads \begin{displaymath}-|\det(u,v_2)\cdot \det A_1'|- |\det(u,v_1) \cdot \det A_2'|+ |\det B_3|+|\det B_4|=0\end{displaymath} we have that the sum of the $\ev\times j$-multiplicities of the curves through a configuration near the wall stays constant.
\item $( (m^2+nm)(\det(u,v_1)+\det(u,v_2))-n\det(v_1,v_2)\big)\leq0$ --- then $\det A_4'$ is positive. A curve of type $\alpha_4$ exists if and only if $\det \tilde{A_i}$ is negative. So in this case $\alpha_1$, $\alpha_2$ and $\alpha_4$ are on one side of the ''wall'' and $\alpha_3$ on the other. But equation \ref{eq1} from above reads \begin{displaymath}-|\det(u,v_2)\cdot \det A_1'|- |\det(u,v_1) \cdot \det A_2'|+ |\det B_3|-|\det B_4|=0\end{displaymath} and we have again that the sum of the $\ev\times j$-multiplicities of the curves through a configuration near the wall stays constant.
\end{itemize}
Let us now come to case \ref{case3}.
As before we can argue that only those curves count, where the $5$-valent vertex is adjacent to the flat cycle. Then the following curves contain $\alpha$ in the boundary and do not count $0$:
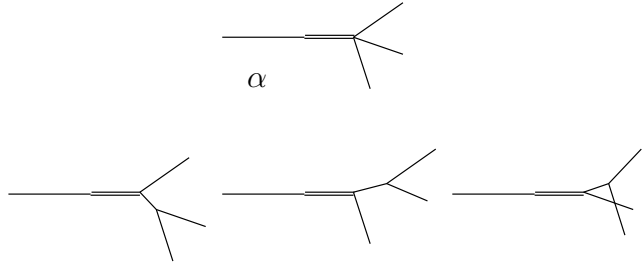
\begin{center}
\input{Graphics/5val.pstex_t}
\end{center}
The proof is here again analogous to case \ref{case1}, only using the ``small'' matrices of $\ev \times j$.

In case \ref{case4}, all curves which have $\alpha$ in their boundary count $0$.
In case \ref{case5}, there is only one possibility with curves that do not count $0$: those where the cycle is adjacent to the $5$-valent vertex. Then the curves which have $\alpha$ in their boundary are the curves where the $4$-valent vertex is resolved, as in \ref{case1}. The proof is analogous to case \ref{case1}, except that we use the ``small'' matrices for $\ev \times j$.

In case \ref{case6}, the $6$-valent vertex has to be adjacent to the cycle, because otherwise every curve which has $\alpha$ in its boundary would count $0$.
So we can now assume that there is a $6$-valent vertex, where two of the adjacent edges are of direction $0$ and form a loop.
To resolve this $6$-valent vertex, we can either form a $5$-valent vertex with a loop and a $3$-valent vertex (these curves are contained in strata of top dimension then), or we can resolve it to four $3$-valent vertices. (We cannot form a flat cycle from the given $6$ edges: the contracted edge must be part of the cycle, and it can either be the whole cycle itself, or it forces the cycle to span $\R^2$.) In the second case, two of the four $3$-valent vertices are connected by the contracted edge and therefore mapped to the same image point in $\R^2$.
Now we want to use the statement that the number of rational curves through given points does not depend on the position of the points for our case here (see \cite{GM051}, respectively use the analogous proof as for proposition 4.4 of \cite{GM053}). More precisely, if there is a point configuration through which a curve with a $4$-valent vertex passes, and we disturb the point configuration slightly, then we always get the same number of tropical curves (counted with multiplicity) passing through the new point configuration.

The image of the $6$-valent vertex (and its adjacent edges) in $\R^2$ looks like a $4$-valent vertex.
The types with one $5$-valent and one $3$-valent vertex are mapped to two $3$-valent vertices, and the type with four $3$-valent vertices is mapped to two $3$-valent vertices and a crossing of two line segments. That is, the images of the $6$-valent vertex as well as of all types which contain it in their boundary look like the possible resolutions of a $4$-valent vertex. We know that there are three types which contain a $4$-valent vertex in their boundary, and we only have to check how we can add contracted bounded edges to these $3$ types, and with which multiplicity they are counted.
The following picture shows the seven possible ways to add contracted bounded edges to the three types:
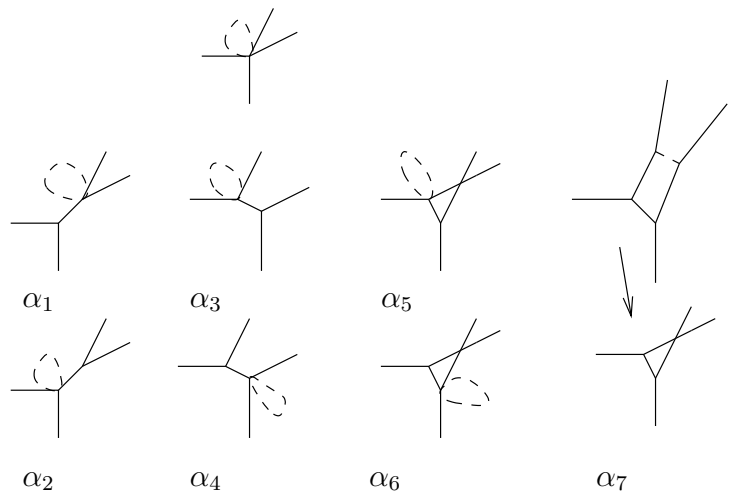
\begin{center}
\input{Graphics/6val.pstex_t}
\end{center}
Note that we can vary the length of the contracted bounded edge in each type, so the curves of these types can have any possible $j$-invariant. The question whether there is a curve of type $\alpha_i$ through a configuration $\mathcal{P}=(p_1,\ldots,p_n,l)\in \R^{2n}\times \overline{\mathcal{M}}_{\trop,\;1,1}$ depends therefore only on the question whether the image of the curve passes through $(p_1,\ldots,p_n)$.
If a curve of type $\alpha_1$ passes through $\mathcal{P}$, then also a curve of type $\alpha_2$ and vice versa. The same holds for $\alpha_3$ and $\alpha_4$, and for $\alpha_5$, $\alpha_6$ and $\alpha_7$. So we only have to see that the sum of the $\ev \times j$-multiplicities of the types whose images are equal can be written as a factor times the multiplicity of the rational curve which arises after removing the contracted bounded edge (and straightening the other edges). Then the statement follows from the statement that the number of rational curves through given points does not depend on the position of the points (see \cite{GM051}, respectively proposition 4.4 of \cite{GM053}).
To see this, we pass to the dual pictures. The $4$-valent vertex is dual to a quadrangle. The images of curves of type $\alpha_1$ and $\alpha_2$ are dual to a subdivision of this quadrangle in two triangles. The same holds for the images of curves of type $\alpha_3$ and $\alpha_4$, however the two triangles arise here by adding the other diagonal. Images of curves of type $\alpha_5$, $\alpha_6$ and $\alpha_7$ are dual to a subdivision consisting of one parallelogram and two triangles.
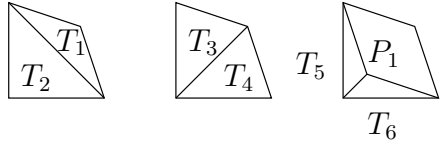
\begin{center}
\input{Graphics/duals.pstex_t}
\end{center}
Lemma \ref{lem-evjcontr-1} tells us that the $\ev \times
j$-multiplicity  of a curve of type $\alpha_1$ is equal to
$(\Area(T_1)-\frac{1}{2})$ times the multiplicity of the rational
curve which arises after removing the contracted bounded edge.
(Recall that the multiplicity of a vertex is by definition equal to
$2\cdot \Area(T)$, where $T$ denotes the dual triangle.) Analogously
the $\ev \times j$-multiplicity of a curve of type $\alpha_2$ is
equal to $(\Area(T_2)-\frac{1}{2})$ times the multiplicity of the
same rational curve. The sum is equal to $(\Area(Q)-1)$ times the
multiplicity of the rational curve, where $Q$ denotes the
quadrangle. We get the same for curves of type $\alpha_3$ and
$\alpha_4$. The sum of the $\ev \times j$-multiplicities of
$\alpha_5$, $\alpha_6$ and $\alpha_7$ is again by lemma
\ref{lem-evjcontr-1} and lemma \ref{lem-evjcontr-2} equal to
$(\Area(T_5)-\frac{1}{2})+(\Area(T_6)-\frac{1}{2})+\Area(P_1)=(\Area(Q)-1)$
times the multiplicity of the corresponding rational curve. Hence
the statement follows.
\end{proof}
\section{Curves with a very large $j$-invariant}
Now we want to use the independence of $\deg_{\ev\times j}(\mathcal{P})$ from $\mathcal{P}$ to compute $\deg_{\ev\times j}(\mathcal{P})$ with the aid of a special configuration $\mathcal{P}=(p_1,\ldots,p_n,l)$ - a configuration where the $j$-invariant $l$ is very large.
\begin {proposition} \label{prop-jbounded}
  Let $n=3d-1$ and $ \mathcal{P}=(p_1,\ldots,p_n,l) \in \R^{2n}\times \overline{\mathcal{M}}_{\trop,\;1,1} $
  be a point in $ \ev\times j $-general position whose $j$-invariant is very
  large (that is, whose image $j(C)\in\overline{\mathcal{M}}_{\trop,\;1,1}$ is a curve with a bounded edge of a very large length). Then every tropical curve $ C \in (\ev\times j)^{-1}(\mathcal{P}) $ with $ \mult_{\ev\times j}(C) \neq 0 $ has
  a contracted bounded edge.
\end {proposition}
\begin {proof}
The proof is similar to proposition 5.1 of \cite{GM053}.
We have to show that the set of all points $ j(C) \in \overline{\mathcal{M}}_{\trop,\;1,1} $ is
  bounded in $ \overline{\mathcal{M}}_{\trop,\;1,1} $, where $C$ runs over all curves in $\mathcal{M}'_{\trop,\;1,n}(d)  $
  with non-zero $ \ev\times j $-multiplicity that have no contracted bounded edge and
  satisfy the given incidence conditions at the marked points. As there are
  only finitely many combinatorial types (analogously to 2.10 of \cite{GM053}) we can
  restrict ourselves to curves of a fixed (but arbitrary) combinatorial type
  $ \alpha $. Since $ \mathcal{P} $ is in $ \ev\times j $-general position we can assume that
  the curves are $3$-valent, respectively contain a flat cycle adjacent to a $4$-valent vertex (it cannot contain a contracted cycle, as we assume that is contains no contracted bounded edges at all).

Assume first that $C$ is $3$-valent. As $C$ is marked by $3d-1$ points we can conclude with 3.50 of \cite{Ma06} that $C$ has a string.
  Analogously to the proof of theorem \ref{thm-degconst} above, we get that there is precisely one string.

  So let $ \Gamma' $ be the unique string in $C$.
Assume first that $\Gamma'\approx \R$.
Then analogously to the proof of 5.1 of \cite{GM053}, we can see that the movement of the string is bounded by the adjacent bounded edges, except if the string consists of only two neighboring unbounded edges. But in this case the only length which is not bounded cannot contribute to the $j$-invariant. So in any case $j(C)$ is bounded.
Assume now $\Gamma'\approx S^1$. As this is the only string, there have to be bounded edges adjacent to the cycle. These bounded edges restrict the movement of the cycle, too. Again $j(C)$ is bounded.
Now assume $C$ has a flat cycle and a $4$-valent vertex adjacent to it, and assume no marked point is adjacent to that $4$-valent vertex. Then all marked points are adjacent to a $3$-valent vertex, and hence we can analogously to 3.50 of \cite{Ma06} see that the curve contains a string. As above, there is exactly one string and its movement is bounded.
Now assume that there is a marked point adjacent to the flat cycle. Then we cannot use 3.50 of \cite{Ma06} to conclude that $C$ has a string. However, the image of the cycle can still not grow arbitrary large:
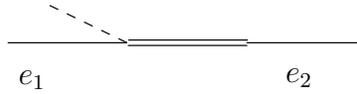
\begin{center}
\input{Graphics/flloop.pstex_t}
\end{center}
The edge $e_2$ has to be bounded: its direction is not a primitive integer vector, as it is equal to the sum of the directions of the two edges of the flat cycle, and therefore it cannot be an unbounded edge. But then the cycle cannot grow arbitrary large.
\end {proof}
As we know that the number of curves $C\in (\ev\times
j)^{-1}(\mathcal{P})$ (counted with multiplicity) does not depend on
$\mathcal{P}$ by \ref{thm-degconst}, we can now choose a special
configuration $\mathcal{P}=(p_1,\ldots,p_n,l)$ where the
$j$-invariant $l$ is very large. Then by \ref{prop-jbounded} we can
conclude that all curves $C\in (\ev\times j)^{-1}(\mathcal{P})$
contain a contracted bounded edge, which is contained in the cycle.
As in lemma \ref{lem-evjcontr-1} and lemma \ref{lem-evjcontr-2},
this contracted bounded edge can either be a loop itself (which is
then adjacent to a $5$-valent vertex), or the contracted bounded
edge is adjacent to two $3$-valent vertices. In both cases, we know
that we can form a rational curve $C'$ of $C$ by removing the
contracted edge and straightening other edges, and we can compute
the $\ev\times j$-multiplicity in terms of the multiplicity of this
rational curve. Note that a rational curve which appears like this
is $3$-valent, as we take an elliptic curve of codimension $0$. The
following lemma shows that we can also ``go back'': we can form
elliptic curves out of a given rational curve which passes through
$(p_1,\ldots,p_n)$.

Define \emph{tropical general position} of the points $(p_1,\ldots,p_n)$ as in \cite{Mi03}, definition 4.7.
Then only simple tropical curves pass through $(p_1,\ldots,p_n)$. In particular, there are only triangles and parallelograms in the dual Newton configurations of these curves.
\begin{lemma}\label{lem-ratel}
Let $n=3d-1$. Take the configuration $\mathcal{P}=(p_1,\ldots,p_n,l)$ where the $j$-invariant $l$ is very large, and such that $(p_1,\ldots,p_n)$ are in tropical general position. Let $C'$ a rational curve which passes through the points $(p_1,\ldots,p_n)$. Then there are several ways to built an elliptic curve $C$ with $j$-invariant $l$ out of $C'$, and the sum of the $\ev\times j$-multiplicities of these elliptic curves is equal to $\binom{d-1}{2}\cdot \mult(C')$.
\end{lemma}
\begin{proof}
Let $V$ be a $3$-valent vertex of $C'$. Then we can make an elliptic curve $C$ out of $C'$ by adding a contracted loop at $V$. There is only one possibility for the length of this loop, as we want to reach that $j(C)=l$. The $\ev \times j$-multiplicity of $C$ is by \ref{lem-evjcontr-1} equal to $\frac{1}{2}(\mult(V)-1)\cdot \mult(C')= (\Area(T)-\frac{1}{2})\cdot \mult(C')$, where $T$ denotes the triangle dual to $V$.
Assume that there is a crossing of two edges $e_1$ and $e_2$ of $C'$, that is, the images $h(e_1)$ and $h(e_2)$ intersect in one point. Then we can add a contracted bounded edge and split $e_1$ and $e_2$ into two edges each. The length of these edges are uniquely determined by the image. The length of the new contracted edge is uniquely determined by the prescribed $j$-invariant $l$. The multiplicity of the elliptic curve $C$ we built like this is due to lemma \ref{lem-evjcontr-2} equal to $|\det(u,v)|\cdot \mult(C')=\Area(P)\cdot \mult(C')$, where $u$ and $v$ denote the two directions of $e_1$ and $e_2$ and $P$ denotes the parallelogram dual to the crossing of $e_1$ and $e_2$.
So it remains to show that $\sum_T (\Area(T)-\frac{1}{2})+ \sum_P \Area(P)=\binom{d-1}{2}$, where $T$ goes over all triangles in the dual Newton subdivision of $C'$ and $P$ goes over all parallelograms.
To see this, we use the theorem of Pick (see \cite{Ful98}, section 5.3).
Denote the number of interior lattice points of a polygon $Q$ by $i(Q)$ and the number of lattice points on the boundary which are not vertices by $b(Q)$.
 Then Pick's formula states that $\Area(T)=i(T)+\frac{b(T)}{2}+\frac{1}{2}$ for a lattice triangle $T$ and $\Area(P)=i(P)+\frac{b(P)}{2}+1$ for a parallelogram $P$.
So we can rewrite the sum from above as
\begin{align*}
&\sum_T \big(\Area(T)-\frac{1}{2}\big)+ \sum_P \Area(P) =\\
&\sum_T \big(i(T)+\frac{b(T)}{2}+\frac{1}{2}-\frac{1}{2}\big)+ \sum_P \big(i(P)+\frac{b(P)}{2}+1\big)=\\
&\sum_T \big(i(T)+\frac{b(T)}{2}\big)+\sum_P \big(i(P)+\frac{b(P)}{2}\big)+\#\{P|P\mbox{ parallelogram in the subdiv}\}=\\
&\sum_T \big(i(T)+\frac{b(T)}{2}\big)+\sum_P \big(i(P)+\frac{b(P)}{2}\big)+\#\{\mbox{lattice points of the subdiv}\}
\end{align*}
where the last equality holds, because $C$ is rational and the genus of a simple tropical curve is equal to the number of points of the subdivision minus the number of parallelograms.
Now we know that the interior lattice points of the big triangle $\Delta_d$ (which is the Newton polygon of curves of degree $d$) that are not contained in the subdivision must either be interior points of a triangle or a parallelogram or on the boundary of a triangle or parallelogram. In the first case, they are counted in  $i(T)$ respectively $i(P)$ of a polygon. In the latter case, as they are interior points of $\Delta_d$, they are part of the boundary of exactly two polygons. That is, in our above sum, they are counted as $b(T)/2$ respectively $b(P)/2$ for two polygons. Hence the first part of the sum counts all interior points which are not part of the subdivision. So we have
\begin{align*}
&\sum_T \big(i(T)+\frac{b(T)}{2}\big)+\sum_P \big(i(P)+\frac{b(P)}{2}\big)+\#\{\mbox{lattice points of the subdiv}\}=\\
&\#\{\mbox{lattice points not contained in the subdiv}\}+ \\& \#\{\mbox{lattice points of the subdiv}\}\\& =
\#\{\mbox{interior points of }\Delta_d\}=\binom{d-1}{2}.
\end{align*}
\end{proof}
We can now sum up our results to prove our main theorem:
\begin{theorem}\label{th-largej}
The number $E_{\trop}(d)$ of tropical elliptic curves passing through $3d-1$ points and with a fixed $j$-invariant, counted with $\ev\times j$-multiplicity, is equal to $\binom{d-1}{2}\cdot N_{\trop}(d)$, where $N_{\trop}(d)$ denotes the number of rational curves through $3d-1$ points (counted with multiplicity).
\end{theorem}
\begin{proof}
The number of tropical elliptic curves passing through $3d-1$ points and with a fixed $j$-invariant is equal to $\deg_{\ev\times j}(\mathcal{P})$, where we can choose any general configuration $\mathcal{P}=(p_1,\ldots,p_n,l)$ due to \ref{thm-degconst}. We choose a configuration with a very large length $l$ as in \ref{prop-jbounded}, and conclude that every elliptic curve passing through this configuration has a contracted bounded edge. From each such elliptic curve with a contracted bounded edge we can form a rational curve by removing the contracted edge and straightening divalent vertices, if necessary. Also, \ref{lem-ratel} tells us that we can go ``backwards'' and form an elliptic curve with $j$-invariant $l$ from each rational curve through $(p_1,\ldots,p_n)$, and that each rational curve contributes with the factor $\binom{d-1}{2}$ to our sum of elliptic curves. Altogether, we have $\binom{d-1}{2}\cdot N_{\trop} (d)$ elliptic curves with $j$-invariant $l$ through $(p_1,\ldots,p_n)$.
\end{proof}
\begin{corollary}
The numbers $E_{\trop}(d)$ and $E(d,j)$ coincide, if $j\notin\{ 0, 1728\}$.
\end{corollary}
\begin{proof}
Theorem \ref{th-largej} tells us that $E_{\trop}(d)=\binom{d-1}{2}N_{\trop}(d)$. The latter is equal to $\binom{d-1}{2}N(d)$ by G.\ Mikhalkin's Correspondence Theorem (see theorem 1 of \cite{Mi03}) and due to Pandharipande's count (\cite{Pan97}), this is equal to $E(d,j)$.
\end{proof}
\section{curves with a very small $j$-invariant}
In the last section, we interpreted a set of rational curves through a given point configuration as elliptic curves with a very large $j$-invariant (and a contracted bounded edge).
Now, we want to interpret the same set of rational curves as elliptic curves with $j$-invariant $0$. Section \ref{sec-mult} helps us to express the multiplicity with which we count the elliptic curves in terms of the rational curves we started with. We will see that we have to count these rational curves with completely different factors than in the previous chapter.
\begin{lemma}\label{lem-ratsmallj}
The number of elliptic curves with a fixed (very small) $j$-invariant and passing through $3d-1$ points in tropical general position is equal to
\begin{displaymath}
\sum_C  \left( \sum_T \big(2\Area(T)^2- \frac{1}{2}\big)\cdot \mult C\right)
\end{displaymath}
where $C$ goes over all rational curves through the $3d-1$ points and $T$ goes over all triangles in the Newton subdivision dual to $C$.
\end{lemma}
\begin{proof}
Given a rational curve $C$, how can we interpret it as an elliptic curve with $j$-invariant $0$? To answer this question, we have to determine how the elliptic curves $C'$ with a very small $j$-invariant which have $C$ in their boundary look like.
There are three possibilities for $C'$:

 Let $\defi(C')=0$. Then the cycle has to disappear to a ($3$-valent) vertex of $C$, hence it must be formed by three edges. Due to lemma \ref{lem-evjcontr-4} the $\ev\times j$-multiplicity is $\mult V\cdot \mult C$ then.
For each vertex $V$ of the rational curve $C$, there are $i(\Delta(V))$ possibilities that a non degenerate cycle disappears to $V$, where $\Delta(V)$ denotes the triangle dual to $V$ and $i(\Delta(V))$ the number of interior points of this triangle. Hence, to count the elliptic curves with a non degenerate cycle we have to count each rational curve $C$ with the factor $\sum_{T} i(T)\cdot 2\Area(T)$ where $T$ goes over all triangles in the Newton subdivision dual to $C$.

Let $\defi(C')=1$. If $e$ is an edge of $C$ with weight bigger $1$, then there can be a small flat cycle at both sides of $e$. The edge $e$ is dual to an edge with interior points in the dual Newton subdivision, and it is in the boundary of two triangles $T_1$ and $T_2$, dual to the two end vertices $V_1$ and $V_2$ of $e$. Assume the flat cycle is adjacent to the vertex $V_1$, and assume that it is formed by two edges with directions $n\cdot u$ and $m\cdot u$, with $\gcd(n,m)=1$ and $(n+m)\cdot u=v$, where $v$ denotes the direction of $e$.
Then by \ref{lem-evjcontr-3} the $\ev\times j$-multiplicity of this curve is \begin{align*} &(n+m)\cdot \det(u,v_1)\cdot \mult C=\det((n+m)u,v_1)\cdot \mult C\\ =& \det(v,v_1)\cdot \mult C=2\Area(T_1)\cdot \mult C,\end{align*} where $v_1$ denotes the direction of another edge adjacent to $V_1$. Respectively, if $n=m=1$ it is \begin{align*} &\det(u,v_1)\cdot \mult C=\det(\frac{1}{2}v,v_1)\cdot \mult C\\&=\Area(T_1)\cdot \mult C.\end{align*}
Assume $\omega(e)$ is even. Then there are $\frac{\omega(e)}{2}-1$ possibilities to separate $e$ to two edges with different directions (that is, with $n\neq m$). Each counts with the factor $2\Area(T_1)$. Also, there is one possibility to separate it to two edges with the same direction, which counts $\Area(T_1)$.
Altogether, we have to count the rational curve with the factor \begin{displaymath}(\frac{\omega(e)}{2}-1)\cdot2\Area(T_1)+\Area(T_1)= (\omega(e)-1)\cdot \Area(T_1).\end{displaymath}
Assume $\omega(e)$ is odd. Then there are $\frac{\omega(e)-1}{2}$ possibilities to split $e$ to two edges with different direction, and each counts with the factor $2\Area(T_1)$.
In any case, we have to count with the factor $(\omega(e)-1)\cdot \Area(T_1)$.
Note that $\omega(e)-1$ is equal to the number of lattice points on the side of the boundary of $T_1$ which is dual to $e$.
But as we have to count these possibilities for all edges of higher weight, we have to add it for all three sides of $T_1$, that is, altogether, we get $b(T_1)\cdot \Area(T_1)$.
Hence, to count the elliptic curves with a flat cycle we have to count each rational curve $C$ with the factor $\sum_{T} b(T)\cdot \Area(T)$ where $T$ goes over all triangles in the Newton subdivision dual to $C$.

Let $\defi(C')=2$. By \ref{lem-evjcontr-1} we have to count it with the factor $\Area(T)-\frac{1}{2}$.
Hence, to count the elliptic curves with a contracted cycle we have to count each rational curve $C$ with the factor $\sum_{T}\Area(T)-\frac{1}{2}$ where $T$ goes over all triangles in the Newton subdivision dual to $C$.

Let us sum up:
\begin{align*}
&\left(\sum_T i(T)\cdot 2\Area(T)+\sum_T b(T)\cdot \Area(T)+\sum_{T}\big(\Area(T)-\frac{1}{2}\big)\right)\cdot \mult C
\\= & \left(\sum_T\big( (2i(T)+b(T)+1) \cdot \Area(T)- \frac{1}{2}\big)\right)\cdot \mult C \\=&
\left( \sum_T\big( 2\Area(T)^2-\frac{1}{2}\big)\right)\cdot \mult C
\end{align*}
where $T$ goes over all triangles in the Newton subdivision dual to $C$.
\end{proof}
At last, we want to apply lemma \ref{lem-ratsmallj} to a set of rational curves passing through a certain point configuration, namely the point configuration which is used in \cite{Mi03}, theorem 2, to prove that marked tropical curves are dual to lattice paths.
Our application results in a faster way to count lattice paths dual to rational curves.

Let $\lambda(x,y)=x-\varepsilon y$ with a very small $\varepsilon >0$. Due to \cite{Mi03} theorem 2 we know that the number of $\lambda$-increasing paths in the triangle $\Delta_d$ is equal to the number of tropical curves through a certain point configuration $\mathcal{P}_{\lambda}$.

Use the notations of chapter 3 of \cite{GM052}.
The tropical curves through $\mathcal{P}_{\lambda}$ are dual to a set of Newton subdivisions. In proposition 3.8 and remark 3.9 of \cite{GM052} we have seen that we can count instead of these Newton subdivisions the column-wise Newton subdivisions for a path. Note that the set of Newton subdivisions which really appear as dual subdivisions of a tropical curve through $\mathcal{P}_{\lambda}$ and the column-wise Newton subdivisions only differ in the location of some parallelograms, the size and locations of the triangles coincide. As for our sum from \ref{lem-ratsmallj} we only count the triangles, we can therefore use the column-wise Newton subdivisions as well.
In remark 3.7 of \cite{GM052} we have seen that a path can only have steps which move one column to the right (with a simultaneous up or down movement), or steps which stay in the same column and move down. The following picture shows such a path and recalls the notations.
\begin{center}
\input{Graphics/not.pstex_t}
\end{center}
For the path in the picture, we have $\alpha^0=5$, $\alpha^1=(1,1) $, $\alpha^2=1 $, $\alpha^3=1$, $\alpha^4=2$, $\alpha^5=0$ and $\alpha^6=0 $; and $h(1)=4 $,  $h(2)=3 $, $h(3)= 2$, $h(4)= 2$, $h(5)=0 $.
The only possibilities for the sequences $\beta'$ are: $\beta'^1=1 $, $\beta'^2= 1$, $\beta'^3= 1$, $\beta'^4=0 $, $\beta'^5=1 $.
The only possibilities for the sequences $\beta$ are: $\beta^0=1 $, $\beta^1=1 $, $\beta^2= 2$, $\beta^3=1 $, $\beta^4=0 $, $\beta^5=0 $.

Proposition 3.8 of \cite{GM052} then gives us a formula to compute the number of column-wise Newton subdivisions times the multiplicity for a path. To get the number we want, we only have to multiply with the factor $( 2\Area(T)^2-\frac{1}{2})$ for each triangle. But note that as in remark 3.9 of \cite{GM052} the position of the triangles below a path are such that they lie in one column and point to the left. That is, they do not have any interior lattice points, and their area is equal to $\frac{1}{2}$ times the length of their right side:
\begin{center}
\input{Graphics/triangle.pstex_t}
\end{center}
There is an analogous statement for triangles above the path, of course.
So, including this factor, we get the following formula: 

\begin{corollary}\label{cor-formrat}
The following formula holds for all $d\geq3$:
\begin{align*}
N(d)&= \frac{1}{\binom{d-1}{2}}\cdot \sum_{\gamma} \sum_{(\beta^0,\ldots,\beta^d),(\beta'^0,\ldots,\beta'^d)}
\binom{\alpha^{i+1}+\beta^{i+1}}{\beta^i}\cdot \binom{\alpha^i+\beta'^i}{\beta'^{i+1}}\\ &\cdot I^{\alpha^{i+1}+\beta^{i+1}-\beta^i}\cdot I^{\alpha^i+\beta'^i-\beta'^{i+1}}\\ &\cdot \left (\frac{I^2-1}{2}\cdot (\alpha^{i+1}+\beta^{i+1}-\beta^i)+\frac{I^2-1}{2}\cdot (\alpha^i+\beta'^i-\beta'^{i+1})\right)
\end{align*}
where the first sum goes over all paths $\gamma$ and the second sum goes over all sequences $(\beta^0,\ldots,\beta^d)$ and $(\beta'^0,\ldots,\beta'^d)$ such that $\beta^0=(d-\alpha^0,0\ldots,0)$, $I \alpha^i+I\beta^i=h(i)$, $\beta'^0=0$ and $d-i-I\beta'^i=h(i)$, and where for a sequence $\alpha=(\alpha_1,\alpha_2,\alpha_3,\ldots)$ $\frac{I^2-1}{2}\cdot \alpha$ denotes the sum $\frac{2^2-1}{2}\cdot \alpha_2+\frac{3^2-1}{2}\cdot \alpha_3 +\ldots$.

\end{corollary}

\begin{proof}
Using G.\ Mikhalkin's Correspondence Theorem (see theorem 2 of \cite{Mi03}) we conclude that $N(d)=N_{\trop}(d)$. Furthermore, $\binom{d-1}{2} N_{{\trop}}(d)=E_{{\trop}}(d)=\deg_{\ev\times j}(\mathcal{P})$, where we can choose any point configuration $\mathcal{P}$ by theorem \ref{th-largej} and theorem \ref{thm-degconst}. 
So $N(d)=\frac{1}{\binom{d-1}{2}}E_{\trop}(d)$ and it remains to argue why the right hand side of the formula above (times $\binom{d-1}{2}$) is equal to $E_{\trop}(d)$. 
We can choose a point $\mathcal{P}=(p_1,\ldots,p_n,l)$ with a very small last coordinate $l$ for the cycle length, and such that $(p_1,\ldots,p_n)$ are in the position described in \cite{Mi03}, theorem 2. We apply lemma \ref{lem-ratsmallj} that tells us that $E_{\trop}(d)=
\sum_C  \left( \sum_T \big(2\Area(T)^2- \frac{1}{2}\big)\cdot \mult C\right)$,
where $C$ goes over all rational curves through the $3d-1$ points and $T$ goes over all triangles in the Newton subdivision dual to $C$.
The Newton subdivision dual to the rational curves through $(p_1,\ldots,p_n)$ differ from the column-wise Newton subdivisions (as defined in remark 3.9 of \cite{GM052}) only in the location of some parallelograms. Size and location of the triangles coincide.
Therefore the above sum is equal to 
$\sum_{N}  \left( \sum_T \big(2\Area(T)^2- \frac{1}{2}\big)\cdot \mult (N)\right)$,
where $N$ goes over all column-wise Newton subdivisions arising from Newton subdivisions dual to rational tropical curves through $(p_1,\ldots,p_n)$.
Proposition 3.8 of \cite{GM052} gives us a formula to compute the number of column-wise Newton subdivisions times their multiplicity. We only have to multiply this formula with the factor $( 2\Area(T)^2-\frac{1}{2})$ for each triangle. As in remark 3.9 of \cite{GM052} the position of the triangles in a column-wise Newton subdivision are such that they lie in one column and point to the left. That is, they do not have any interior lattice points, and their area is equal to $\frac{1}{2}$ times the length of their right side.
The factor $\binom{\alpha^{i+1}+\beta^{i+1}}{\beta^i}\cdot \binom{\alpha^i+\beta'^i}{\beta'^{i+1}}$ counts the possibilities to arrange parallelograms below and above the path (hence the number of Newton subdivisions). The factor $I^{\alpha^{i+1}+\beta^{i+1}-\beta^i}\cdot I^{\alpha^i+\beta'^i-\beta'^{i+1}}$ counts the double areas of the triangles - hence the multiplicity of the curves dual to the path. (See also remark 3.9 of \cite{GM052}).
The factor $\frac{I^2-1}{2}\cdot (\alpha^{i+1}+\beta^{i+1}-\beta^i)+\frac{I^2-1}{2}\cdot (\alpha^i+\beta'^i-\beta'^{i+1})$ is the factor $( 2\Area(T)^2-\frac{1}{2})$ for each triangle.
\end{proof}

Note that even though this sum looks at the first glance more complicated than the sum from proposition 3.8 of \cite{GM052}, it is easier to compute, because we count a lot of paths with the factor $0$ --- all paths with only steps of size $1$.

\begin{example}
For $d=3$, there is only one lattice path with a step of size bigger
than one.
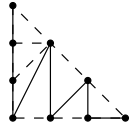
\begin{center}
\input{Graphics/d3.pstex_t}
\end{center}
There is only one possible Newton subdivision for this path, as indicated in the picture. There are two triangles of area $1$. Both contribute $\frac{3}{2} \cdot \mult C=\frac{3}{2} \cdot 4=6$. Altogether, we get $6+6=12=N(3)$, as expected.
\end{example}
\begin{example}
For $d=4$, we only have to consider the paths below, because all other paths have only steps of size $1$.
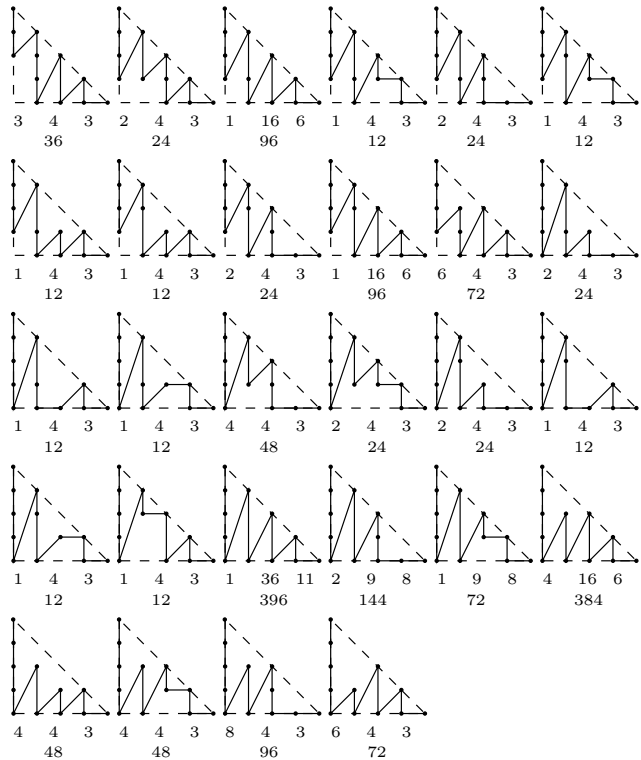
\begin{center}
\input{Graphics/d4.pstex_t}
\end{center}
There are three numbers in the first row below each path: the first number is the number of possible Newton subdivisions. The second number is the multiplicity of the tropical curves dual to these Newton subdivisions. (Hence the product of the first two numbers is the multiplicity of the path.) The third number is the factor $\sum_T\big( 2\Area(T)^2-\frac{1}{2}\big)$ with which we have to count here.
The fourth number, in the second row, is the product of the three numbers above, so we have to count each path with that number.
The sum of the numbers is the second row is $1860=3\cdot 620=3\cdot N(d)$, as claimed.
\end{example}

\end{document}

%% file: Graphics/absg1.pstex_t
\begin{picture}(0,0)%
\includegraphics{Graphics/absg1.pstex}%
\end{picture}%
\setlength{\unitlength}{3947sp}%
\begingroup\makeatletter\ifx\SetFigFont\undefined%
\gdef\SetFigFont#1#2#3#4#5{%
  \reset@font\fontsize{#1}{#2pt}%
  \fontfamily{#3}\fontseries{#4}\fontshape{#5}%
  \selectfont}%
\fi\endgroup%
\begin{picture}(1824,324)(4339,-3373)
\end{picture}%

%% file: Graphics/boundary3.pstex_t
\begin{picture}(0,0)%
\includegraphics{Graphics/boundary3.pstex}%
\end{picture}%
\setlength{\unitlength}{3947sp}%
\begingroup\makeatletter\ifx\SetFigFont\undefined%
\gdef\SetFigFont#1#2#3#4#5{%
  \reset@font\fontsize{#1}{#2pt}%
  \fontfamily{#3}\fontseries{#4}\fontshape{#5}%
  \selectfont}%
\fi\endgroup%
\begin{picture}(4883,1142)(5659,-4388)
\put(7348,-4188){\makebox(0,0)[lb]{\smash{{\SetFigFont{8}{9.6}{\familydefault}{\mddefault}{\updefault}{\color[rgb]{0,0,0}$\Gamma_1$}%
}}}}
\put(8762,-4356){\makebox(0,0)[lb]{\smash{{\SetFigFont{8}{9.6}{\familydefault}{\mddefault}{\updefault}{\color[rgb]{0,0,0}$\Gamma$}%
}}}}
\put(6016,-4299){\makebox(0,0)[lb]{\smash{{\SetFigFont{8}{9.6}{\familydefault}{\mddefault}{\updefault}{\color[rgb]{0,0,0}$\Gamma$}%
}}}}
\put(9977,-4274){\makebox(0,0)[lb]{\smash{{\SetFigFont{8}{9.6}{\familydefault}{\mddefault}{\updefault}{\color[rgb]{0,0,0}$\Gamma_1$}%
}}}}
\end{picture}%

%% file: Graphics/flatloop.pstex_t
\begin{picture}(0,0)%
\includegraphics{Graphics/flatloop.pstex}%
\end{picture}%
\setlength{\unitlength}{3947sp}%
\begingroup\makeatletter\ifx\SetFigFont\undefined%
\gdef\SetFigFont#1#2#3#4#5{%
  \reset@font\fontsize{#1}{#2pt}%
  \fontfamily{#3}\fontseries{#4}\fontshape{#5}%
  \selectfont}%
\fi\endgroup%
\begin{picture}(2124,1343)(3889,-4242)
\put(5476,-3061){\makebox(0,0)[lb]{\smash{{\SetFigFont{12}{14.4}{\familydefault}{\mddefault}{\updefault}{\color[rgb]{0,0,0}$v$}%
}}}}
\put(4812,-3326){\makebox(0,0)[lb]{\smash{{\SetFigFont{12}{14.4}{\familydefault}{\mddefault}{\updefault}{\color[rgb]{0,0,0}$n\cdot u$}%
}}}}
\put(4853,-4193){\makebox(0,0)[lb]{\smash{{\SetFigFont{12}{14.4}{\familydefault}{\mddefault}{\updefault}{\color[rgb]{0,0,0}$m\cdot u$}%
}}}}
\end{picture}%

%% file: Graphics/five.pstex_t
\begin{picture}(0,0)%
\includegraphics{Graphics/five.pstex}%
\end{picture}%
\setlength{\unitlength}{3947sp}%
\begingroup\makeatletter\ifx\SetFigFont\undefined%
\gdef\SetFigFont#1#2#3#4#5{%
  \reset@font\fontsize{#1}{#2pt}%
  \fontfamily{#3}\fontseries{#4}\fontshape{#5}%
  \selectfont}%
\fi\endgroup%
\begin{picture}(1910,1794)(4639,-5173)
\put(5026,-3811){\makebox(0,0)[lb]{\smash{{\SetFigFont{12}{14.4}{\familydefault}{\mddefault}{\updefault}{\color[rgb]{0,0,0}$u$}%
}}}}
\put(5701,-3511){\makebox(0,0)[lb]{\smash{{\SetFigFont{12}{14.4}{\familydefault}{\mddefault}{\updefault}{\color[rgb]{0,0,0}$v$}%
}}}}
\end{picture}%

%% file: Graphics/elr.pstex_t
\begin{picture}(0,0)%
\includegraphics{Graphics/elr.pstex}%
\end{picture}%
\setlength{\unitlength}{3947sp}%
\begingroup\makeatletter\ifx\SetFigFont\undefined%
\gdef\SetFigFont#1#2#3#4#5{%
  \reset@font\fontsize{#1}{#2pt}%
  \fontfamily{#3}\fontseries{#4}\fontshape{#5}%
  \selectfont}%
\fi\endgroup%
\begin{picture}(3624,1374)(2239,-4423)
\put(3908,-3669){\makebox(0,0)[lb]{\smash{{\SetFigFont{12}{14.4}{\familydefault}{\mddefault}{\updefault}{\color[rgb]{0,0,0}$h$}%
}}}}
\end{picture}%

%% file: Graphics/elr1.pstex_t
\begin{picture}(0,0)%
\includegraphics{Graphics/elr1.pstex}%
\end{picture}%
\setlength{\unitlength}{3947sp}%
\begingroup\makeatletter\ifx\SetFigFont\undefined%
\gdef\SetFigFont#1#2#3#4#5{%
  \reset@font\fontsize{#1}{#2pt}%
  \fontfamily{#3}\fontseries{#4}\fontshape{#5}%
  \selectfont}%
\fi\endgroup%
\begin{picture}(1957,964)(2089,-3987)
\end{picture}%

%% file: Graphics/jex.pstex_t
\begin{picture}(0,0)%
\includegraphics{Graphics/jex.pstex}%
\end{picture}%
\setlength{\unitlength}{3947sp}%
\begingroup\makeatletter\ifx\SetFigFont\undefined%
\gdef\SetFigFont#1#2#3#4#5{%
  \reset@font\fontsize{#1}{#2pt}%
  \fontfamily{#3}\fontseries{#4}\fontshape{#5}%
  \selectfont}%
\fi\endgroup%
\begin{picture}(5385,1892)(11803,-4948)
\put(12117,-3182){\makebox(0,0)[lb]{\smash{{\SetFigFont{10}{12.0}{\familydefault}{\mddefault}{\updefault}{\color[rgb]{0,0,0}$x_1$}%
}}}}
\put(12843,-4029){\makebox(0,0)[lb]{\smash{{\SetFigFont{10}{12.0}{\familydefault}{\mddefault}{\updefault}{\color[rgb]{0,0,0}$l_1$}%
}}}}
\put(13266,-4029){\makebox(0,0)[lb]{\smash{{\SetFigFont{10}{12.0}{\familydefault}{\mddefault}{\updefault}{\color[rgb]{0,0,0}$l_2$}%
}}}}
\put(13085,-4755){\makebox(0,0)[lb]{\smash{{\SetFigFont{10}{12.0}{\familydefault}{\mddefault}{\updefault}{\color[rgb]{0,0,0}$l_5$}%
}}}}
\put(12783,-4694){\makebox(0,0)[lb]{\smash{{\SetFigFont{10}{12.0}{\familydefault}{\mddefault}{\updefault}{\color[rgb]{0,0,0}$l_6$}%
}}}}
\put(12601,-4513){\makebox(0,0)[lb]{\smash{{\SetFigFont{10}{12.0}{\familydefault}{\mddefault}{\updefault}{\color[rgb]{0,0,0}$l_7$}%
}}}}
\put(13569,-4271){\makebox(0,0)[lb]{\smash{{\SetFigFont{10}{12.0}{\familydefault}{\mddefault}{\updefault}{\color[rgb]{0,0,0}$l_3$}%
}}}}
\put(12601,-4271){\makebox(0,0)[lb]{\smash{{\SetFigFont{10}{12.0}{\familydefault}{\mddefault}{\updefault}{\color[rgb]{0,0,0}$l_8$}%
}}}}
\put(13448,-4573){\makebox(0,0)[lb]{\smash{{\SetFigFont{10}{12.0}{\familydefault}{\mddefault}{\updefault}{\color[rgb]{0,0,0}$l_4$}%
}}}}
\put(15976,-3811){\makebox(0,0)[lb]{\smash{{\SetFigFont{12}{14.4}{\familydefault}{\mddefault}{\updefault}{\color[rgb]{0,0,0}$j(C)$}%
}}}}
\end{picture}%

%% file: Graphics/evj.pstex_t
\begin{picture}(0,0)%
\includegraphics{Graphics/evj.pstex}%
\end{picture}%
\setlength{\unitlength}{3947sp}%
\begingroup\makeatletter\ifx\SetFigFont\undefined%
\gdef\SetFigFont#1#2#3#4#5{%
  \reset@font\fontsize{#1}{#2pt}%
  \fontfamily{#3}\fontseries{#4}\fontshape{#5}%
  \selectfont}%
\fi\endgroup%
\begin{picture}(2332,1271)(3915,-4624)
\put(4633,-4414){\makebox(0,0)[lb]{\smash{{\SetFigFont{8}{9.6}{\familydefault}{\mddefault}{\updefault}{\color[rgb]{0,0,0}$V$}%
}}}}
\put(5406,-3863){\makebox(0,0)[lb]{\smash{{\SetFigFont{8}{9.6}{\familydefault}{\mddefault}{\updefault}{\color[rgb]{0,0,0}$v_4$}%
}}}}
\put(5129,-4526){\makebox(0,0)[lb]{\smash{{\SetFigFont{8}{9.6}{\familydefault}{\mddefault}{\updefault}{\color[rgb]{0,0,0}$v_5$}%
}}}}
\put(4412,-4582){\makebox(0,0)[lb]{\smash{{\SetFigFont{8}{9.6}{\familydefault}{\mddefault}{\updefault}{\color[rgb]{0,0,0}$v_1$}%
}}}}
\put(5461,-3476){\makebox(0,0)[lb]{\smash{{\SetFigFont{8}{9.6}{\familydefault}{\mddefault}{\updefault}{\color[rgb]{0,0,0}$x_2$}%
}}}}
\put(3915,-4029){\makebox(0,0)[lb]{\smash{{\SetFigFont{8}{9.6}{\familydefault}{\mddefault}{\updefault}{\color[rgb]{0,0,0}$x_1$}%
}}}}
\put(4579,-3917){\makebox(0,0)[lb]{\smash{{\SetFigFont{8}{9.6}{\familydefault}{\mddefault}{\updefault}{\color[rgb]{0,0,0}$v_2$}%
}}}}
\put(4855,-3641){\makebox(0,0)[lb]{\smash{{\SetFigFont{8}{9.6}{\familydefault}{\mddefault}{\updefault}{\color[rgb]{0,0,0}$v_3$}%
}}}}
\end{picture}%

%% file: Graphics/contr.pstex_t
\begin{picture}(0,0)%
\includegraphics{Graphics/contr.pstex}%
\end{picture}%
\setlength{\unitlength}{3947sp}%
\begingroup\makeatletter\ifx\SetFigFont\undefined%
\gdef\SetFigFont#1#2#3#4#5{%
  \reset@font\fontsize{#1}{#2pt}%
  \fontfamily{#3}\fontseries{#4}\fontshape{#5}%
  \selectfont}%
\fi\endgroup%
\begin{picture}(4188,1344)(4189,-3823)
\put(6001,-3061){\makebox(0,0)[lb]{\smash{{\SetFigFont{12}{14.4}{\familydefault}{\mddefault}{\updefault}{\color[rgb]{0,0,0}$h$}%
}}}}
\put(4876,-3286){\makebox(0,0)[lb]{\smash{{\SetFigFont{12}{14.4}{\familydefault}{\mddefault}{\updefault}{\color[rgb]{0,0,0}$e$}%
}}}}
\put(4876,-2836){\makebox(0,0)[lb]{\smash{{\SetFigFont{12}{14.4}{\familydefault}{\mddefault}{\updefault}{\color[rgb]{0,0,0}$e_1$}%
}}}}
\put(5101,-3586){\makebox(0,0)[lb]{\smash{{\SetFigFont{12}{14.4}{\familydefault}{\mddefault}{\updefault}{\color[rgb]{0,0,0}$e_4$}%
}}}}
\put(4276,-3736){\makebox(0,0)[lb]{\smash{{\SetFigFont{12}{14.4}{\familydefault}{\mddefault}{\updefault}{\color[rgb]{0,0,0}$e_3$}%
}}}}
\put(4276,-2986){\makebox(0,0)[lb]{\smash{{\SetFigFont{12}{14.4}{\familydefault}{\mddefault}{\updefault}{\color[rgb]{0,0,0}$e_2$}%
}}}}
\put(7801,-2686){\makebox(0,0)[lb]{\smash{{\SetFigFont{12}{14.4}{\familydefault}{\mddefault}{\updefault}{\color[rgb]{0,0,0}$u$}%
}}}}
\put(6901,-3211){\makebox(0,0)[lb]{\smash{{\SetFigFont{12}{14.4}{\familydefault}{\mddefault}{\updefault}{\color[rgb]{0,0,0}$u$}%
}}}}
\put(7201,-2611){\makebox(0,0)[lb]{\smash{{\SetFigFont{12}{14.4}{\familydefault}{\mddefault}{\updefault}{\color[rgb]{0,0,0}$v$}%
}}}}
\put(8101,-3211){\makebox(0,0)[lb]{\smash{{\SetFigFont{12}{14.4}{\familydefault}{\mddefault}{\updefault}{\color[rgb]{0,0,0}$v$}%
}}}}
\end{picture}%

%% file: Graphics/fl.pstex_t
\begin{picture}(0,0)%
\includegraphics{Graphics/fl.pstex}%
\end{picture}%
\setlength{\unitlength}{3947sp}%
\begingroup\makeatletter\ifx\SetFigFont\undefined%
\gdef\SetFigFont#1#2#3#4#5{%
  \reset@font\fontsize{#1}{#2pt}%
  \fontfamily{#3}\fontseries{#4}\fontshape{#5}%
  \selectfont}%
\fi\endgroup%
\begin{picture}(1824,1231)(3814,-3560)
\put(3901,-2461){\makebox(0,0)[lb]{\smash{{\SetFigFont{12}{14.4}{\familydefault}{\mddefault}{\updefault}{\color[rgb]{0,0,0}$v$}%
}}}}
\put(4276,-2686){\makebox(0,0)[lb]{\smash{{\SetFigFont{12}{14.4}{\familydefault}{\mddefault}{\updefault}{\color[rgb]{0,0,0}$m\cdot u$}%
}}}}
\put(4876,-2686){\makebox(0,0)[lb]{\smash{{\SetFigFont{12}{14.4}{\familydefault}{\mddefault}{\updefault}{\color[rgb]{0,0,0}$(m+n)\cdot u$}%
}}}}
\put(4351,-3511){\makebox(0,0)[lb]{\smash{{\SetFigFont{12}{14.4}{\familydefault}{\mddefault}{\updefault}{\color[rgb]{0,0,0}$n\cdot u$}%
}}}}
\end{picture}%

%% file: Graphics/dreieck.pstex_t
\begin{picture}(0,0)%
\includegraphics{Graphics/dreieck.pstex}%
\end{picture}%
\setlength{\unitlength}{3947sp}%
\begingroup\makeatletter\ifx\SetFigFont\undefined%
\gdef\SetFigFont#1#2#3#4#5{%
  \reset@font\fontsize{#1}{#2pt}%
  \fontfamily{#3}\fontseries{#4}\fontshape{#5}%
  \selectfont}%
\fi\endgroup%
\begin{picture}(2885,1636)(5089,-4685)
\put(5401,-4636){\makebox(0,0)[lb]{\smash{{\SetFigFont{12}{14.4}{\familydefault}{\mddefault}{\updefault}{\color[rgb]{0,0,0}$C$}%
}}}}
\put(7501,-4636){\makebox(0,0)[lb]{\smash{{\SetFigFont{12}{14.4}{\familydefault}{\mddefault}{\updefault}{\color[rgb]{0,0,0}$C'$}%
}}}}
\put(5926,-3586){\makebox(0,0)[lb]{\smash{{\SetFigFont{12}{14.4}{\familydefault}{\mddefault}{\updefault}{\color[rgb]{0,0,0}$v_1$}%
}}}}
\put(5476,-4111){\makebox(0,0)[lb]{\smash{{\SetFigFont{12}{14.4}{\familydefault}{\mddefault}{\updefault}{\color[rgb]{0,0,0}$v_3$}%
}}}}
\put(7651,-3886){\makebox(0,0)[lb]{\smash{{\SetFigFont{12}{14.4}{\familydefault}{\mddefault}{\updefault}{\color[rgb]{0,0,0}$V$}%
}}}}
\put(5176,-3436){\makebox(0,0)[lb]{\smash{{\SetFigFont{12}{14.4}{\familydefault}{\mddefault}{\updefault}{\color[rgb]{0,0,0}$v_2$}%
}}}}
\end{picture}%

%% file: Graphics/dickesdreieck.pstex_t
\begin{picture}(0,0)%
\includegraphics{Graphics/dickesdreieck.pstex}%
\end{picture}%
\setlength{\unitlength}{3947sp}%
\begingroup\makeatletter\ifx\SetFigFont\undefined%
\gdef\SetFigFont#1#2#3#4#5{%
  \reset@font\fontsize{#1}{#2pt}%
  \fontfamily{#3}\fontseries{#4}\fontshape{#5}%
  \selectfont}%
\fi\endgroup%
\begin{picture}(2724,624)(5239,-5923)
\end{picture}%

%% file: Graphics/four.pstex_t
\begin{picture}(0,0)%
\includegraphics{Graphics/four.pstex}%
\end{picture}%
\setlength{\unitlength}{3947sp}%
\begingroup\makeatletter\ifx\SetFigFont\undefined%
\gdef\SetFigFont#1#2#3#4#5{%
  \reset@font\fontsize{#1}{#2pt}%
  \fontfamily{#3}\fontseries{#4}\fontshape{#5}%
  \selectfont}%
\fi\endgroup%
\begin{picture}(3909,1102)(5989,-4001)
\put(6264,-3961){\makebox(0,0)[lb]{\smash{{\SetFigFont{8}{9.6}{\familydefault}{\mddefault}{\updefault}{\color[rgb]{0,0,0}$\alpha$}%
}}}}
\put(7366,-3961){\makebox(0,0)[lb]{\smash{{\SetFigFont{8}{9.6}{\familydefault}{\mddefault}{\updefault}{\color[rgb]{0,0,0}$\alpha_1$}%
}}}}
\put(8469,-3961){\makebox(0,0)[lb]{\smash{{\SetFigFont{8}{9.6}{\familydefault}{\mddefault}{\updefault}{\color[rgb]{0,0,0}$\alpha_2$}%
}}}}
\put(9624,-3961){\makebox(0,0)[lb]{\smash{{\SetFigFont{8}{9.6}{\familydefault}{\mddefault}{\updefault}{\color[rgb]{0,0,0}$\alpha_3$}%
}}}}
\put(7471,-3383){\makebox(0,0)[lb]{\smash{{\SetFigFont{8}{9.6}{\familydefault}{\mddefault}{\updefault}{\color[rgb]{0,0,0}$e$}%
}}}}
\put(8469,-3331){\makebox(0,0)[lb]{\smash{{\SetFigFont{8}{9.6}{\familydefault}{\mddefault}{\updefault}{\color[rgb]{0,0,0}$e$}%
}}}}
\put(9519,-3488){\makebox(0,0)[lb]{\smash{{\SetFigFont{8}{9.6}{\familydefault}{\mddefault}{\updefault}{\color[rgb]{0,0,0}$e$}%
}}}}
\end{picture}%

%% file: Graphics/string14.pstex_t
\begin{picture}(0,0)%
\includegraphics{Graphics/string14.pstex}%
\end{picture}%
\setlength{\unitlength}{3947sp}%
\begingroup\makeatletter\ifx\SetFigFont\undefined%
\gdef\SetFigFont#1#2#3#4#5{%
  \reset@font\fontsize{#1}{#2pt}%
  \fontfamily{#3}\fontseries{#4}\fontshape{#5}%
  \selectfont}%
\fi\endgroup%
\begin{picture}(3500,2920)(4313,-5744)
\put(4876,-4711){\makebox(0,0)[lb]{\smash{{\SetFigFont{12}{14.4}{\familydefault}{\mddefault}{\updefault}{\color[rgb]{0,0,0}$(3)$}%
}}}}
\put(6676,-4711){\makebox(0,0)[lb]{\smash{{\SetFigFont{12}{14.4}{\familydefault}{\mddefault}{\updefault}{\color[rgb]{0,0,0}$(4)$}%
}}}}
\put(4876,-3661){\makebox(0,0)[lb]{\smash{{\SetFigFont{12}{14.4}{\familydefault}{\mddefault}{\updefault}{\color[rgb]{0,0,0}$(1)$}%
}}}}
\put(6676,-3661){\makebox(0,0)[lb]{\smash{{\SetFigFont{12}{14.4}{\familydefault}{\mddefault}{\updefault}{\color[rgb]{0,0,0}$(2)$}%
}}}}
\put(4876,-5686){\makebox(0,0)[lb]{\smash{{\SetFigFont{12}{14.4}{\familydefault}{\mddefault}{\updefault}{\color[rgb]{0,0,0}$(5)$}%
}}}}
\put(6676,-5686){\makebox(0,0)[lb]{\smash{{\SetFigFont{12}{14.4}{\familydefault}{\mddefault}{\updefault}{\color[rgb]{0,0,0}$(6)$}%
}}}}
\end{picture}%

%% file: Graphics/4types1.pstex_t
\begin{picture}(0,0)%
\includegraphics{Graphics/4types1.pstex}%
\end{picture}%
\setlength{\unitlength}{3947sp}%
\begingroup\makeatletter\ifx\SetFigFont\undefined%
\gdef\SetFigFont#1#2#3#4#5{%
  \reset@font\fontsize{#1}{#2pt}%
  \fontfamily{#3}\fontseries{#4}\fontshape{#5}%
  \selectfont}%
\fi\endgroup%
\begin{picture}(2799,1489)(5164,-1694)
\put(6376,-361){\makebox(0,0)[lb]{\smash{{\SetFigFont{12}{14.4}{\familydefault}{\mddefault}{\updefault}{\color[rgb]{0,0,0}$x_1$}%
}}}}
\put(6526,-961){\makebox(0,0)[lb]{\smash{{\SetFigFont{12}{14.4}{\familydefault}{\mddefault}{\updefault}{\color[rgb]{0,0,0}$e_6$}%
}}}}
\put(6001,-961){\makebox(0,0)[lb]{\smash{{\SetFigFont{12}{14.4}{\familydefault}{\mddefault}{\updefault}{\color[rgb]{0,0,0}$e_5$}%
}}}}
\put(6301,-1336){\makebox(0,0)[lb]{\smash{{\SetFigFont{12}{14.4}{\familydefault}{\mddefault}{\updefault}{\color[rgb]{0,0,0}$e_7$}%
}}}}
\put(5401,-661){\makebox(0,0)[lb]{\smash{{\SetFigFont{12}{14.4}{\familydefault}{\mddefault}{\updefault}{\color[rgb]{0,0,0}$e_1$}%
}}}}
\put(5401,-1636){\makebox(0,0)[lb]{\smash{{\SetFigFont{12}{14.4}{\familydefault}{\mddefault}{\updefault}{\color[rgb]{0,0,0}$e_4$}%
}}}}
\put(7276,-1636){\makebox(0,0)[lb]{\smash{{\SetFigFont{12}{14.4}{\familydefault}{\mddefault}{\updefault}{\color[rgb]{0,0,0}$e_3$}%
}}}}
\put(7276,-736){\makebox(0,0)[lb]{\smash{{\SetFigFont{12}{14.4}{\familydefault}{\mddefault}{\updefault}{\color[rgb]{0,0,0}$e_2$}%
}}}}
\end{picture}%

%% file: Graphics/4types.pstex_t
\begin{picture}(0,0)%
\includegraphics{Graphics/4types.pstex}%
\end{picture}%
\setlength{\unitlength}{3947sp}%
\begingroup\makeatletter\ifx\SetFigFont\undefined%
\gdef\SetFigFont#1#2#3#4#5{%
  \reset@font\fontsize{#1}{#2pt}%
  \fontfamily{#3}\fontseries{#4}\fontshape{#5}%
  \selectfont}%
\fi\endgroup%
\begin{picture}(3653,2646)(4684,-4694)
\put(5026,-2386){\makebox(0,0)[lb]{\smash{{\SetFigFont{12}{14.4}{\familydefault}{\mddefault}{\updefault}{\color[rgb]{0,0,0}$V$}%
}}}}
\put(6826,-2686){\makebox(0,0)[lb]{\smash{{\SetFigFont{12}{14.4}{\familydefault}{\mddefault}{\updefault}{\color[rgb]{0,0,0}$V$}%
}}}}
\put(5026,-3586){\makebox(0,0)[lb]{\smash{{\SetFigFont{12}{14.4}{\familydefault}{\mddefault}{\updefault}{\color[rgb]{0,0,0}$V$}%
}}}}
\put(6826,-3586){\makebox(0,0)[lb]{\smash{{\SetFigFont{12}{14.4}{\familydefault}{\mddefault}{\updefault}{\color[rgb]{0,0,0}$V$}%
}}}}
\put(6676,-3886){\makebox(0,0)[lb]{\smash{{\SetFigFont{12}{14.4}{\familydefault}{\mddefault}{\updefault}{\color[rgb]{0,0,0}$e_8$}%
}}}}
\put(7876,-3886){\makebox(0,0)[lb]{\smash{{\SetFigFont{12}{14.4}{\familydefault}{\mddefault}{\updefault}{\color[rgb]{0,0,0}$e_9$}%
}}}}
\put(6076,-3886){\makebox(0,0)[lb]{\smash{{\SetFigFont{12}{14.4}{\familydefault}{\mddefault}{\updefault}{\color[rgb]{0,0,0}$e_9$}%
}}}}
\put(4801,-3886){\makebox(0,0)[lb]{\smash{{\SetFigFont{12}{14.4}{\familydefault}{\mddefault}{\updefault}{\color[rgb]{0,0,0}$e_8$}%
}}}}
\put(5326,-2986){\makebox(0,0)[lb]{\smash{{\SetFigFont{12}{14.4}{\familydefault}{\mddefault}{\updefault}{\color[rgb]{0,0,0}$\alpha_1$}%
}}}}
\put(7051,-2986){\makebox(0,0)[lb]{\smash{{\SetFigFont{12}{14.4}{\familydefault}{\mddefault}{\updefault}{\color[rgb]{0,0,0}$\alpha_2$}%
}}}}
\put(5176,-4561){\makebox(0,0)[lb]{\smash{{\SetFigFont{12}{14.4}{\familydefault}{\mddefault}{\updefault}{\color[rgb]{0,0,0}$\alpha_3$}%
}}}}
\put(7126,-4636){\makebox(0,0)[lb]{\smash{{\SetFigFont{12}{14.4}{\familydefault}{\mddefault}{\updefault}{\color[rgb]{0,0,0}$\alpha_4$}%
}}}}
\put(7551,-2373){\makebox(0,0)[lb]{\smash{{\SetFigFont{12}{14.4}{\familydefault}{\mddefault}{\updefault}{\color[rgb]{0,0,0}$e_8$}%
}}}}
\put(5061,-2677){\makebox(0,0)[lb]{\smash{{\SetFigFont{12}{14.4}{\familydefault}{\mddefault}{\updefault}{\color[rgb]{0,0,0}$e_8$}%
}}}}
\end{picture}%

%% file: Graphics/5val.pstex_t
\begin{picture}(0,0)%
\includegraphics{Graphics/5val.pstex}%
\end{picture}%
\setlength{\unitlength}{3947sp}%
\begingroup\makeatletter\ifx\SetFigFont\undefined%
\gdef\SetFigFont#1#2#3#4#5{%
  \reset@font\fontsize{#1}{#2pt}%
  \fontfamily{#3}\fontseries{#4}\fontshape{#5}%
  \selectfont}%
\fi\endgroup%
\begin{picture}(3999,1648)(5539,-3947)
\put(7051,-2836){\makebox(0,0)[lb]{\smash{{\SetFigFont{12}{14.4}{\familydefault}{\mddefault}{\updefault}{\color[rgb]{0,0,0}$\alpha$}%
}}}}
\end{picture}%

%% file: Graphics/6val.pstex_t
\begin{picture}(0,0)%
\includegraphics{Graphics/6val.pstex}%
\end{picture}%
\setlength{\unitlength}{3947sp}%
\begingroup\makeatletter\ifx\SetFigFont\undefined%
\gdef\SetFigFont#1#2#3#4#5{%
  \reset@font\fontsize{#1}{#2pt}%
  \fontfamily{#3}\fontseries{#4}\fontshape{#5}%
  \selectfont}%
\fi\endgroup%
\begin{picture}(4524,3070)(2989,-5819)
\put(3076,-4636){\makebox(0,0)[lb]{\smash{{\SetFigFont{12}{14.4}{\familydefault}{\mddefault}{\updefault}{\color[rgb]{0,0,0}$\alpha_1$}%
}}}}
\put(3076,-5761){\makebox(0,0)[lb]{\smash{{\SetFigFont{12}{14.4}{\familydefault}{\mddefault}{\updefault}{\color[rgb]{0,0,0}$\alpha_2$}%
}}}}
\put(4126,-4636){\makebox(0,0)[lb]{\smash{{\SetFigFont{12}{14.4}{\familydefault}{\mddefault}{\updefault}{\color[rgb]{0,0,0}$\alpha_3$}%
}}}}
\put(4126,-5761){\makebox(0,0)[lb]{\smash{{\SetFigFont{12}{14.4}{\familydefault}{\mddefault}{\updefault}{\color[rgb]{0,0,0}$\alpha_4$}%
}}}}
\put(5326,-4636){\makebox(0,0)[lb]{\smash{{\SetFigFont{12}{14.4}{\familydefault}{\mddefault}{\updefault}{\color[rgb]{0,0,0}$\alpha_5$}%
}}}}
\put(5251,-5761){\makebox(0,0)[lb]{\smash{{\SetFigFont{12}{14.4}{\familydefault}{\mddefault}{\updefault}{\color[rgb]{0,0,0}$\alpha_6$}%
}}}}
\put(6676,-5761){\makebox(0,0)[lb]{\smash{{\SetFigFont{12}{14.4}{\familydefault}{\mddefault}{\updefault}{\color[rgb]{0,0,0}$\alpha_7$}%
}}}}
\end{picture}%

%% file: Graphics/duals.pstex_t
\begin{picture}(0,0)%
\includegraphics{Graphics/duals.pstex}%
\end{picture}%
\setlength{\unitlength}{3947sp}%
\begingroup\makeatletter\ifx\SetFigFont\undefined%
\gdef\SetFigFont#1#2#3#4#5{%
  \reset@font\fontsize{#1}{#2pt}%
  \fontfamily{#3}\fontseries{#4}\fontshape{#5}%
  \selectfont}%
\fi\endgroup%
\begin{picture}(2738,895)(5239,-3944)
\put(7501,-3436){\makebox(0,0)[lb]{\smash{{\SetFigFont{12}{14.4}{\familydefault}{\mddefault}{\updefault}{\color[rgb]{0,0,0}$P_1$}%
}}}}
\put(7051,-3511){\makebox(0,0)[lb]{\smash{{\SetFigFont{12}{14.4}{\familydefault}{\mddefault}{\updefault}{\color[rgb]{0,0,0}$T_5$}%
}}}}
\put(7501,-3886){\makebox(0,0)[lb]{\smash{{\SetFigFont{12}{14.4}{\familydefault}{\mddefault}{\updefault}{\color[rgb]{0,0,0}$T_6$}%
}}}}
\put(5551,-3361){\makebox(0,0)[lb]{\smash{{\SetFigFont{12}{14.4}{\familydefault}{\mddefault}{\updefault}{\color[rgb]{0,0,0}$T_1$}%
}}}}
\put(5326,-3586){\makebox(0,0)[lb]{\smash{{\SetFigFont{12}{14.4}{\familydefault}{\mddefault}{\updefault}{\color[rgb]{0,0,0}$T_2$}%
}}}}
\put(6376,-3361){\makebox(0,0)[lb]{\smash{{\SetFigFont{12}{14.4}{\familydefault}{\mddefault}{\updefault}{\color[rgb]{0,0,0}$T_3$}%
}}}}
\put(6601,-3586){\makebox(0,0)[lb]{\smash{{\SetFigFont{12}{14.4}{\familydefault}{\mddefault}{\updefault}{\color[rgb]{0,0,0}$T_4$}%
}}}}
\end{picture}%

%% file: Graphics/flloop.pstex_t
\begin{picture}(0,0)%
\includegraphics{Graphics/flloop.pstex}%
\end{picture}%
\setlength{\unitlength}{3947sp}%
\begingroup\makeatletter\ifx\SetFigFont\undefined%
\gdef\SetFigFont#1#2#3#4#5{%
  \reset@font\fontsize{#1}{#2pt}%
  \fontfamily{#3}\fontseries{#4}\fontshape{#5}%
  \selectfont}%
\fi\endgroup%
\begin{picture}(2274,577)(5689,-3827)
\put(7446,-3761){\makebox(0,0)[lb]{\smash{{\SetFigFont{12}{14.4}{\familydefault}{\mddefault}{\updefault}{\color[rgb]{0,0,0}$e_2$}%
}}}}
\put(5769,-3769){\makebox(0,0)[lb]{\smash{{\SetFigFont{12}{14.4}{\familydefault}{\mddefault}{\updefault}{\color[rgb]{0,0,0}$e_1$}%
}}}}
\end{picture}%

%% file: Graphics/not.pstex_t
\begin{picture}(0,0)%
\includegraphics{Graphics/not.pstex}%
\end{picture}%
\setlength{\unitlength}{3947sp}%
\begingroup\makeatletter\ifx\SetFigFont\undefined%
\gdef\SetFigFont#1#2#3#4#5{%
  \reset@font\fontsize{#1}{#2pt}%
  \fontfamily{#3}\fontseries{#4}\fontshape{#5}%
  \selectfont}%
\fi\endgroup%
\begin{picture}(1918,1874)(4126,-4148)
\put(5491,-3276){\makebox(0,0)[lb]{\smash{{\SetFigFont{7}{8.4}{\familydefault}{\mddefault}{\updefault}{\color[rgb]{0,0,0}$h(4)$}%
}}}}
\put(5751,-3991){\makebox(0,0)[lb]{\smash{{\SetFigFont{7}{8.4}{\familydefault}{\mddefault}{\updefault}{\color[rgb]{0,0,0}$h(5)$}%
}}}}
\put(4256,-3081){\makebox(0,0)[lb]{\smash{{\SetFigFont{7}{8.4}{\familydefault}{\mddefault}{\updefault}{\color[rgb]{0,0,0}$\alpha^0$}%
}}}}
\put(5296,-3731){\makebox(0,0)[lb]{\smash{{\SetFigFont{7}{8.4}{\familydefault}{\mddefault}{\updefault}{\color[rgb]{0,0,0}$\alpha^4$}%
}}}}
\put(4451,-4121){\makebox(0,0)[lb]{\smash{{\SetFigFont{7}{8.4}{\familydefault}{\mddefault}{\updefault}{\color[rgb]{0,0,0}$\beta^1$}%
}}}}
\put(4126,-4121){\makebox(0,0)[lb]{\smash{{\SetFigFont{7}{8.4}{\familydefault}{\mddefault}{\updefault}{\color[rgb]{0,0,0}$\beta^0$}%
}}}}
\put(5101,-4121){\makebox(0,0)[lb]{\smash{{\SetFigFont{7}{8.4}{\familydefault}{\mddefault}{\updefault}{\color[rgb]{0,0,0}$\beta^2$}%
}}}}
\put(4906,-2496){\makebox(0,0)[lb]{\smash{{\SetFigFont{7}{8.4}{\familydefault}{\mddefault}{\updefault}{\color[rgb]{0,0,0}$\beta'^1$}%
}}}}
\put(5036,-3146){\makebox(0,0)[lb]{\smash{{\SetFigFont{7}{8.4}{\familydefault}{\mddefault}{\updefault}{\color[rgb]{0,0,0}$h(2)$}%
}}}}
\put(5101,-3276){\makebox(0,0)[lb]{\smash{{\SetFigFont{7}{8.4}{\familydefault}{\mddefault}{\updefault}{\color[rgb]{0,0,0}$h(3)$}%
}}}}
\put(5166,-2821){\makebox(0,0)[lb]{\smash{{\SetFigFont{7}{8.4}{\familydefault}{\mddefault}{\updefault}{\color[rgb]{0,0,0}$\beta'^2$}%
}}}}
\put(4728,-2909){\makebox(0,0)[lb]{\smash{{\SetFigFont{7}{8.4}{\familydefault}{\mddefault}{\updefault}{\color[rgb]{0,0,0}$h(1)$}%
}}}}
\put(4529,-3464){\makebox(0,0)[lb]{\smash{{\SetFigFont{7}{8.4}{\familydefault}{\mddefault}{\updefault}{\color[rgb]{0,0,0}$\alpha^1$}%
}}}}
\end{picture}%

%% file: Graphics/triangle.pstex_t
\begin{picture}(0,0)%
\includegraphics{Graphics/triangle.pstex}%
\end{picture}%
\setlength{\unitlength}{3947sp}%
\begingroup\makeatletter\ifx\SetFigFont\undefined%
\gdef\SetFigFont#1#2#3#4#5{%
  \reset@font\fontsize{#1}{#2pt}%
  \fontfamily{#3}\fontseries{#4}\fontshape{#5}%
  \selectfont}%
\fi\endgroup%
\begin{picture}(774,774)(4939,-4723)
\end{picture}%

%% file: Graphics/d3.pstex_t
\begin{picture}(0,0)%
\includegraphics{Graphics/d3.pstex}%
\end{picture}%
\setlength{\unitlength}{3947sp}%
\begingroup\makeatletter\ifx\SetFigFont\undefined%
\gdef\SetFigFont#1#2#3#4#5{%
  \reset@font\fontsize{#1}{#2pt}%
  \fontfamily{#3}\fontseries{#4}\fontshape{#5}%
  \selectfont}%
\fi\endgroup%
\begin{picture}(761,761)(4008,-3404)
\end{picture}%

%% file: Graphics/d4.pstex_t
\begin{picture}(0,0)%
\includegraphics{Graphics/d4.pstex}%
\end{picture}%
\setlength{\unitlength}{3947sp}%
\begingroup\makeatletter\ifx\SetFigFont\undefined%
\gdef\SetFigFont#1#2#3#4#5{%
  \reset@font\fontsize{#1}{#2pt}%
  \fontfamily{#3}\fontseries{#4}\fontshape{#5}%
  \selectfont}%
\fi\endgroup%
\begin{picture}(3950,4747)(3809,-11278)
\put(6487,-7287){\makebox(0,0)[lb]{\smash{{\SetFigFont{6}{7.2}{\familydefault}{\mddefault}{\updefault}{\color[rgb]{0,0,0}$2$}%
}}}}
\put(7151,-7287){\makebox(0,0)[lb]{\smash{{\SetFigFont{6}{7.2}{\familydefault}{\mddefault}{\updefault}{\color[rgb]{0,0,0}$1$}%
}}}}
\put(3829,-7287){\makebox(0,0)[lb]{\smash{{\SetFigFont{6}{7.2}{\familydefault}{\mddefault}{\updefault}{\color[rgb]{0,0,0}$3$}%
}}}}
\put(4271,-7287){\makebox(0,0)[lb]{\smash{{\SetFigFont{6}{7.2}{\familydefault}{\mddefault}{\updefault}{\color[rgb]{0,0,0}$3$}%
}}}}
\put(4493,-7287){\makebox(0,0)[lb]{\smash{{\SetFigFont{6}{7.2}{\familydefault}{\mddefault}{\updefault}{\color[rgb]{0,0,0}$2$}%
}}}}
\put(4715,-7287){\makebox(0,0)[lb]{\smash{{\SetFigFont{6}{7.2}{\familydefault}{\mddefault}{\updefault}{\color[rgb]{0,0,0}$4$}%
}}}}
\put(4936,-7287){\makebox(0,0)[lb]{\smash{{\SetFigFont{6}{7.2}{\familydefault}{\mddefault}{\updefault}{\color[rgb]{0,0,0}$3$}%
}}}}
\put(5158,-7287){\makebox(0,0)[lb]{\smash{{\SetFigFont{6}{7.2}{\familydefault}{\mddefault}{\updefault}{\color[rgb]{0,0,0}$1$}%
}}}}
\put(5821,-7287){\makebox(0,0)[lb]{\smash{{\SetFigFont{6}{7.2}{\familydefault}{\mddefault}{\updefault}{\color[rgb]{0,0,0}$1$}%
}}}}
\put(6044,-7287){\makebox(0,0)[lb]{\smash{{\SetFigFont{6}{7.2}{\familydefault}{\mddefault}{\updefault}{\color[rgb]{0,0,0}$4$}%
}}}}
\put(6265,-7287){\makebox(0,0)[lb]{\smash{{\SetFigFont{6}{7.2}{\familydefault}{\mddefault}{\updefault}{\color[rgb]{0,0,0}$3$}%
}}}}
\put(6708,-7287){\makebox(0,0)[lb]{\smash{{\SetFigFont{6}{7.2}{\familydefault}{\mddefault}{\updefault}{\color[rgb]{0,0,0}$4$}%
}}}}
\put(6929,-7287){\makebox(0,0)[lb]{\smash{{\SetFigFont{6}{7.2}{\familydefault}{\mddefault}{\updefault}{\color[rgb]{0,0,0}$3$}%
}}}}
\put(7372,-7287){\makebox(0,0)[lb]{\smash{{\SetFigFont{6}{7.2}{\familydefault}{\mddefault}{\updefault}{\color[rgb]{0,0,0}$4$}%
}}}}
\put(7594,-7287){\makebox(0,0)[lb]{\smash{{\SetFigFont{6}{7.2}{\familydefault}{\mddefault}{\updefault}{\color[rgb]{0,0,0}$3$}%
}}}}
\put(3829,-8246){\makebox(0,0)[lb]{\smash{{\SetFigFont{6}{7.2}{\familydefault}{\mddefault}{\updefault}{\color[rgb]{0,0,0}$1$}%
}}}}
\put(4051,-8246){\makebox(0,0)[lb]{\smash{{\SetFigFont{6}{7.2}{\familydefault}{\mddefault}{\updefault}{\color[rgb]{0,0,0}$4$}%
}}}}
\put(4271,-8246){\makebox(0,0)[lb]{\smash{{\SetFigFont{6}{7.2}{\familydefault}{\mddefault}{\updefault}{\color[rgb]{0,0,0}$3$}%
}}}}
\put(4493,-8246){\makebox(0,0)[lb]{\smash{{\SetFigFont{6}{7.2}{\familydefault}{\mddefault}{\updefault}{\color[rgb]{0,0,0}$1$}%
}}}}
\put(4715,-8246){\makebox(0,0)[lb]{\smash{{\SetFigFont{6}{7.2}{\familydefault}{\mddefault}{\updefault}{\color[rgb]{0,0,0}$4$}%
}}}}
\put(4936,-8246){\makebox(0,0)[lb]{\smash{{\SetFigFont{6}{7.2}{\familydefault}{\mddefault}{\updefault}{\color[rgb]{0,0,0}$3$}%
}}}}
\put(5158,-8246){\makebox(0,0)[lb]{\smash{{\SetFigFont{6}{7.2}{\familydefault}{\mddefault}{\updefault}{\color[rgb]{0,0,0}$2$}%
}}}}
\put(5380,-8246){\makebox(0,0)[lb]{\smash{{\SetFigFont{6}{7.2}{\familydefault}{\mddefault}{\updefault}{\color[rgb]{0,0,0}$4$}%
}}}}
\put(5601,-8246){\makebox(0,0)[lb]{\smash{{\SetFigFont{6}{7.2}{\familydefault}{\mddefault}{\updefault}{\color[rgb]{0,0,0}$3$}%
}}}}
\put(5821,-8246){\makebox(0,0)[lb]{\smash{{\SetFigFont{6}{7.2}{\familydefault}{\mddefault}{\updefault}{\color[rgb]{0,0,0}$1$}%
}}}}
\put(6044,-8246){\makebox(0,0)[lb]{\smash{{\SetFigFont{6}{7.2}{\familydefault}{\mddefault}{\updefault}{\color[rgb]{0,0,0}$16$}%
}}}}
\put(6265,-8246){\makebox(0,0)[lb]{\smash{{\SetFigFont{6}{7.2}{\familydefault}{\mddefault}{\updefault}{\color[rgb]{0,0,0}$6$}%
}}}}
\put(6708,-8246){\makebox(0,0)[lb]{\smash{{\SetFigFont{6}{7.2}{\familydefault}{\mddefault}{\updefault}{\color[rgb]{0,0,0}$4$}%
}}}}
\put(6929,-8246){\makebox(0,0)[lb]{\smash{{\SetFigFont{6}{7.2}{\familydefault}{\mddefault}{\updefault}{\color[rgb]{0,0,0}$3$}%
}}}}
\put(7151,-8246){\makebox(0,0)[lb]{\smash{{\SetFigFont{6}{7.2}{\familydefault}{\mddefault}{\updefault}{\color[rgb]{0,0,0}$2$}%
}}}}
\put(7372,-8246){\makebox(0,0)[lb]{\smash{{\SetFigFont{6}{7.2}{\familydefault}{\mddefault}{\updefault}{\color[rgb]{0,0,0}$4$}%
}}}}
\put(7594,-8246){\makebox(0,0)[lb]{\smash{{\SetFigFont{6}{7.2}{\familydefault}{\mddefault}{\updefault}{\color[rgb]{0,0,0}$3$}%
}}}}
\put(3829,-9207){\makebox(0,0)[lb]{\smash{{\SetFigFont{6}{7.2}{\familydefault}{\mddefault}{\updefault}{\color[rgb]{0,0,0}$1$}%
}}}}
\put(4051,-9207){\makebox(0,0)[lb]{\smash{{\SetFigFont{6}{7.2}{\familydefault}{\mddefault}{\updefault}{\color[rgb]{0,0,0}$4$}%
}}}}
\put(4271,-9207){\makebox(0,0)[lb]{\smash{{\SetFigFont{6}{7.2}{\familydefault}{\mddefault}{\updefault}{\color[rgb]{0,0,0}$3$}%
}}}}
\put(4493,-9207){\makebox(0,0)[lb]{\smash{{\SetFigFont{6}{7.2}{\familydefault}{\mddefault}{\updefault}{\color[rgb]{0,0,0}$1$}%
}}}}
\put(4715,-9207){\makebox(0,0)[lb]{\smash{{\SetFigFont{6}{7.2}{\familydefault}{\mddefault}{\updefault}{\color[rgb]{0,0,0}$4$}%
}}}}
\put(4936,-9207){\makebox(0,0)[lb]{\smash{{\SetFigFont{6}{7.2}{\familydefault}{\mddefault}{\updefault}{\color[rgb]{0,0,0}$3$}%
}}}}
\put(5158,-9207){\makebox(0,0)[lb]{\smash{{\SetFigFont{6}{7.2}{\familydefault}{\mddefault}{\updefault}{\color[rgb]{0,0,0}$4$}%
}}}}
\put(5380,-9207){\makebox(0,0)[lb]{\smash{{\SetFigFont{6}{7.2}{\familydefault}{\mddefault}{\updefault}{\color[rgb]{0,0,0}$4$}%
}}}}
\put(5601,-9207){\makebox(0,0)[lb]{\smash{{\SetFigFont{6}{7.2}{\familydefault}{\mddefault}{\updefault}{\color[rgb]{0,0,0}$3$}%
}}}}
\put(5821,-9207){\makebox(0,0)[lb]{\smash{{\SetFigFont{6}{7.2}{\familydefault}{\mddefault}{\updefault}{\color[rgb]{0,0,0}$2$}%
}}}}
\put(6044,-9207){\makebox(0,0)[lb]{\smash{{\SetFigFont{6}{7.2}{\familydefault}{\mddefault}{\updefault}{\color[rgb]{0,0,0}$4$}%
}}}}
\put(6265,-9207){\makebox(0,0)[lb]{\smash{{\SetFigFont{6}{7.2}{\familydefault}{\mddefault}{\updefault}{\color[rgb]{0,0,0}$3$}%
}}}}
\put(6487,-9207){\makebox(0,0)[lb]{\smash{{\SetFigFont{6}{7.2}{\familydefault}{\mddefault}{\updefault}{\color[rgb]{0,0,0}$2$}%
}}}}
\put(6708,-9207){\makebox(0,0)[lb]{\smash{{\SetFigFont{6}{7.2}{\familydefault}{\mddefault}{\updefault}{\color[rgb]{0,0,0}$4$}%
}}}}
\put(6929,-9207){\makebox(0,0)[lb]{\smash{{\SetFigFont{6}{7.2}{\familydefault}{\mddefault}{\updefault}{\color[rgb]{0,0,0}$3$}%
}}}}
\put(7151,-9207){\makebox(0,0)[lb]{\smash{{\SetFigFont{6}{7.2}{\familydefault}{\mddefault}{\updefault}{\color[rgb]{0,0,0}$1$}%
}}}}
\put(7372,-9207){\makebox(0,0)[lb]{\smash{{\SetFigFont{6}{7.2}{\familydefault}{\mddefault}{\updefault}{\color[rgb]{0,0,0}$4$}%
}}}}
\put(7594,-9207){\makebox(0,0)[lb]{\smash{{\SetFigFont{6}{7.2}{\familydefault}{\mddefault}{\updefault}{\color[rgb]{0,0,0}$3$}%
}}}}
\put(3829,-10166){\makebox(0,0)[lb]{\smash{{\SetFigFont{6}{7.2}{\familydefault}{\mddefault}{\updefault}{\color[rgb]{0,0,0}$1$}%
}}}}
\put(4051,-10166){\makebox(0,0)[lb]{\smash{{\SetFigFont{6}{7.2}{\familydefault}{\mddefault}{\updefault}{\color[rgb]{0,0,0}$4$}%
}}}}
\put(4271,-10166){\makebox(0,0)[lb]{\smash{{\SetFigFont{6}{7.2}{\familydefault}{\mddefault}{\updefault}{\color[rgb]{0,0,0}$3$}%
}}}}
\put(4493,-10166){\makebox(0,0)[lb]{\smash{{\SetFigFont{6}{7.2}{\familydefault}{\mddefault}{\updefault}{\color[rgb]{0,0,0}$1$}%
}}}}
\put(4715,-10166){\makebox(0,0)[lb]{\smash{{\SetFigFont{6}{7.2}{\familydefault}{\mddefault}{\updefault}{\color[rgb]{0,0,0}$4$}%
}}}}
\put(4936,-10166){\makebox(0,0)[lb]{\smash{{\SetFigFont{6}{7.2}{\familydefault}{\mddefault}{\updefault}{\color[rgb]{0,0,0}$3$}%
}}}}
\put(5158,-10166){\makebox(0,0)[lb]{\smash{{\SetFigFont{6}{7.2}{\familydefault}{\mddefault}{\updefault}{\color[rgb]{0,0,0}$1$}%
}}}}
\put(5380,-10166){\makebox(0,0)[lb]{\smash{{\SetFigFont{6}{7.2}{\familydefault}{\mddefault}{\updefault}{\color[rgb]{0,0,0}$36$}%
}}}}
\put(5601,-10166){\makebox(0,0)[lb]{\smash{{\SetFigFont{6}{7.2}{\familydefault}{\mddefault}{\updefault}{\color[rgb]{0,0,0}$11$}%
}}}}
\put(5821,-10166){\makebox(0,0)[lb]{\smash{{\SetFigFont{6}{7.2}{\familydefault}{\mddefault}{\updefault}{\color[rgb]{0,0,0}$2$}%
}}}}
\put(6044,-10166){\makebox(0,0)[lb]{\smash{{\SetFigFont{6}{7.2}{\familydefault}{\mddefault}{\updefault}{\color[rgb]{0,0,0}$9$}%
}}}}
\put(6265,-10166){\makebox(0,0)[lb]{\smash{{\SetFigFont{6}{7.2}{\familydefault}{\mddefault}{\updefault}{\color[rgb]{0,0,0}$8$}%
}}}}
\put(6487,-10166){\makebox(0,0)[lb]{\smash{{\SetFigFont{6}{7.2}{\familydefault}{\mddefault}{\updefault}{\color[rgb]{0,0,0}$1$}%
}}}}
\put(6708,-10166){\makebox(0,0)[lb]{\smash{{\SetFigFont{6}{7.2}{\familydefault}{\mddefault}{\updefault}{\color[rgb]{0,0,0}$9$}%
}}}}
\put(6929,-10166){\makebox(0,0)[lb]{\smash{{\SetFigFont{6}{7.2}{\familydefault}{\mddefault}{\updefault}{\color[rgb]{0,0,0}$8$}%
}}}}
\put(7151,-10166){\makebox(0,0)[lb]{\smash{{\SetFigFont{6}{7.2}{\familydefault}{\mddefault}{\updefault}{\color[rgb]{0,0,0}$4$}%
}}}}
\put(7372,-10166){\makebox(0,0)[lb]{\smash{{\SetFigFont{6}{7.2}{\familydefault}{\mddefault}{\updefault}{\color[rgb]{0,0,0}$16$}%
}}}}
\put(7594,-10166){\makebox(0,0)[lb]{\smash{{\SetFigFont{6}{7.2}{\familydefault}{\mddefault}{\updefault}{\color[rgb]{0,0,0}$6$}%
}}}}
\put(3829,-11127){\makebox(0,0)[lb]{\smash{{\SetFigFont{6}{7.2}{\familydefault}{\mddefault}{\updefault}{\color[rgb]{0,0,0}$4$}%
}}}}
\put(4051,-11127){\makebox(0,0)[lb]{\smash{{\SetFigFont{6}{7.2}{\familydefault}{\mddefault}{\updefault}{\color[rgb]{0,0,0}$4$}%
}}}}
\put(4271,-11127){\makebox(0,0)[lb]{\smash{{\SetFigFont{6}{7.2}{\familydefault}{\mddefault}{\updefault}{\color[rgb]{0,0,0}$3$}%
}}}}
\put(4493,-11127){\makebox(0,0)[lb]{\smash{{\SetFigFont{6}{7.2}{\familydefault}{\mddefault}{\updefault}{\color[rgb]{0,0,0}$4$}%
}}}}
\put(4715,-11127){\makebox(0,0)[lb]{\smash{{\SetFigFont{6}{7.2}{\familydefault}{\mddefault}{\updefault}{\color[rgb]{0,0,0}$4$}%
}}}}
\put(4936,-11127){\makebox(0,0)[lb]{\smash{{\SetFigFont{6}{7.2}{\familydefault}{\mddefault}{\updefault}{\color[rgb]{0,0,0}$3$}%
}}}}
\put(5158,-11127){\makebox(0,0)[lb]{\smash{{\SetFigFont{6}{7.2}{\familydefault}{\mddefault}{\updefault}{\color[rgb]{0,0,0}$8$}%
}}}}
\put(5380,-11127){\makebox(0,0)[lb]{\smash{{\SetFigFont{6}{7.2}{\familydefault}{\mddefault}{\updefault}{\color[rgb]{0,0,0}$4$}%
}}}}
\put(5601,-11127){\makebox(0,0)[lb]{\smash{{\SetFigFont{6}{7.2}{\familydefault}{\mddefault}{\updefault}{\color[rgb]{0,0,0}$3$}%
}}}}
\put(6044,-11127){\makebox(0,0)[lb]{\smash{{\SetFigFont{6}{7.2}{\familydefault}{\mddefault}{\updefault}{\color[rgb]{0,0,0}$4$}%
}}}}
\put(6265,-11127){\makebox(0,0)[lb]{\smash{{\SetFigFont{6}{7.2}{\familydefault}{\mddefault}{\updefault}{\color[rgb]{0,0,0}$3$}%
}}}}
\put(6487,-8246){\makebox(0,0)[lb]{\smash{{\SetFigFont{6}{7.2}{\familydefault}{\mddefault}{\updefault}{\color[rgb]{0,0,0}$6$}%
}}}}
\put(5380,-7287){\makebox(0,0)[lb]{\smash{{\SetFigFont{6}{7.2}{\familydefault}{\mddefault}{\updefault}{\color[rgb]{0,0,0}$16$}%
}}}}
\put(5601,-7287){\makebox(0,0)[lb]{\smash{{\SetFigFont{6}{7.2}{\familydefault}{\mddefault}{\updefault}{\color[rgb]{0,0,0}$6$}%
}}}}
\put(4051,-7287){\makebox(0,0)[lb]{\smash{{\SetFigFont{6}{7.2}{\familydefault}{\mddefault}{\updefault}{\color[rgb]{0,0,0}$4$}%
}}}}
\put(5821,-11127){\makebox(0,0)[lb]{\smash{{\SetFigFont{6}{7.2}{\familydefault}{\mddefault}{\updefault}{\color[rgb]{0,0,0}$6$}%
}}}}
\put(4018,-7415){\makebox(0,0)[lb]{\smash{{\SetFigFont{6}{7.2}{\familydefault}{\mddefault}{\updefault}{\color[rgb]{0,0,0}$36$}%
}}}}
\put(4695,-7415){\makebox(0,0)[lb]{\smash{{\SetFigFont{6}{7.2}{\familydefault}{\mddefault}{\updefault}{\color[rgb]{0,0,0}$24$}%
}}}}
\put(5373,-7415){\makebox(0,0)[lb]{\smash{{\SetFigFont{6}{7.2}{\familydefault}{\mddefault}{\updefault}{\color[rgb]{0,0,0}$96$}%
}}}}
\put(6050,-7415){\makebox(0,0)[lb]{\smash{{\SetFigFont{6}{7.2}{\familydefault}{\mddefault}{\updefault}{\color[rgb]{0,0,0}$12$}%
}}}}
\put(6671,-7415){\makebox(0,0)[lb]{\smash{{\SetFigFont{6}{7.2}{\familydefault}{\mddefault}{\updefault}{\color[rgb]{0,0,0}$24$}%
}}}}
\put(7350,-7415){\makebox(0,0)[lb]{\smash{{\SetFigFont{6}{7.2}{\familydefault}{\mddefault}{\updefault}{\color[rgb]{0,0,0}$12$}%
}}}}
\put(4018,-8374){\makebox(0,0)[lb]{\smash{{\SetFigFont{6}{7.2}{\familydefault}{\mddefault}{\updefault}{\color[rgb]{0,0,0}$12$}%
}}}}
\put(4695,-8374){\makebox(0,0)[lb]{\smash{{\SetFigFont{6}{7.2}{\familydefault}{\mddefault}{\updefault}{\color[rgb]{0,0,0}$12$}%
}}}}
\put(5373,-8374){\makebox(0,0)[lb]{\smash{{\SetFigFont{6}{7.2}{\familydefault}{\mddefault}{\updefault}{\color[rgb]{0,0,0}$24$}%
}}}}
\put(6050,-8374){\makebox(0,0)[lb]{\smash{{\SetFigFont{6}{7.2}{\familydefault}{\mddefault}{\updefault}{\color[rgb]{0,0,0}$96$}%
}}}}
\put(6671,-8374){\makebox(0,0)[lb]{\smash{{\SetFigFont{6}{7.2}{\familydefault}{\mddefault}{\updefault}{\color[rgb]{0,0,0}$72$}%
}}}}
\put(7350,-8374){\makebox(0,0)[lb]{\smash{{\SetFigFont{6}{7.2}{\familydefault}{\mddefault}{\updefault}{\color[rgb]{0,0,0}$24$}%
}}}}
\put(4018,-9335){\makebox(0,0)[lb]{\smash{{\SetFigFont{6}{7.2}{\familydefault}{\mddefault}{\updefault}{\color[rgb]{0,0,0}$12$}%
}}}}
\put(4695,-9335){\makebox(0,0)[lb]{\smash{{\SetFigFont{6}{7.2}{\familydefault}{\mddefault}{\updefault}{\color[rgb]{0,0,0}$12$}%
}}}}
\put(5373,-9335){\makebox(0,0)[lb]{\smash{{\SetFigFont{6}{7.2}{\familydefault}{\mddefault}{\updefault}{\color[rgb]{0,0,0}$48$}%
}}}}
\put(6050,-9335){\makebox(0,0)[lb]{\smash{{\SetFigFont{6}{7.2}{\familydefault}{\mddefault}{\updefault}{\color[rgb]{0,0,0}$24$}%
}}}}
\put(6728,-9335){\makebox(0,0)[lb]{\smash{{\SetFigFont{6}{7.2}{\familydefault}{\mddefault}{\updefault}{\color[rgb]{0,0,0}$24$}%
}}}}
\put(7350,-9335){\makebox(0,0)[lb]{\smash{{\SetFigFont{6}{7.2}{\familydefault}{\mddefault}{\updefault}{\color[rgb]{0,0,0}$12$}%
}}}}
\put(4018,-10294){\makebox(0,0)[lb]{\smash{{\SetFigFont{6}{7.2}{\familydefault}{\mddefault}{\updefault}{\color[rgb]{0,0,0}$12$}%
}}}}
\put(4695,-10294){\makebox(0,0)[lb]{\smash{{\SetFigFont{6}{7.2}{\familydefault}{\mddefault}{\updefault}{\color[rgb]{0,0,0}$12$}%
}}}}
\put(5373,-10294){\makebox(0,0)[lb]{\smash{{\SetFigFont{6}{7.2}{\familydefault}{\mddefault}{\updefault}{\color[rgb]{0,0,0}$396$}%
}}}}
\put(5994,-10294){\makebox(0,0)[lb]{\smash{{\SetFigFont{6}{7.2}{\familydefault}{\mddefault}{\updefault}{\color[rgb]{0,0,0}$144$}%
}}}}
\put(6671,-10294){\makebox(0,0)[lb]{\smash{{\SetFigFont{6}{7.2}{\familydefault}{\mddefault}{\updefault}{\color[rgb]{0,0,0}$72$}%
}}}}
\put(7350,-10294){\makebox(0,0)[lb]{\smash{{\SetFigFont{6}{7.2}{\familydefault}{\mddefault}{\updefault}{\color[rgb]{0,0,0}$384$}%
}}}}
\put(4018,-11255){\makebox(0,0)[lb]{\smash{{\SetFigFont{6}{7.2}{\familydefault}{\mddefault}{\updefault}{\color[rgb]{0,0,0}$48$}%
}}}}
\put(4695,-11255){\makebox(0,0)[lb]{\smash{{\SetFigFont{6}{7.2}{\familydefault}{\mddefault}{\updefault}{\color[rgb]{0,0,0}$48$}%
}}}}
\put(5373,-11255){\makebox(0,0)[lb]{\smash{{\SetFigFont{6}{7.2}{\familydefault}{\mddefault}{\updefault}{\color[rgb]{0,0,0}$96$}%
}}}}
\put(6050,-11255){\makebox(0,0)[lb]{\smash{{\SetFigFont{6}{7.2}{\familydefault}{\mddefault}{\updefault}{\color[rgb]{0,0,0}$72$}%
}}}}
\end{picture}%